\documentclass[final,10pt]{m2an}

\usepackage{subfig}
\usepackage{graphicx}
\usepackage{amssymb}
\usepackage{amsfonts}
\usepackage{amsmath}
\usepackage{lineno}
\usepackage{bbm}

\newtheorem{ssmptn}[thrm]{Assumption}

\begin{document}

\title{Measuring the Irreversibility of Numerical Schemes for Reversible Stochastic 
Differential Equations} \thanks{We thanks Natesh Pillai for useful comments and suggestions. M. A. K. and Y.P. are partially supported by NSF-CMMI 0835673 and 
L. R.-B. is partially supported by NSF -DMS-1109316}

\author{Markos Katsoulakis}
\address{Department of Mathematics and Statistics, University of Massachusetts, Amherst, MA, USA.}
\author{Yannis Pantazis}\sameaddress{1}
\author{Luc Rey-Bellet} \sameaddress{1}

\begin{abstract} For a stationary Markov process the detailed balance condition is
equivalent to the time-reversibility of the process. For stochastic differential equations (SDE's),
the time discretization of numerical schemes usually destroys the 
time-reversibility property. Despite an extensive literature on the numerical analysis for  SDE's,  their 
stability properties, strong and/or weak error estimates, large deviations and infinite-time 
estimates, no quantitative results are known on the lack of reversibility of discrete-time 
approximation processes. In this paper we provide such quantitative estimates by using 
the concept of entropy production rate, inspired by ideas from non-equilibrium statistical 
mechanics.  The entropy production rate for a stochastic process is defined as the relative 
entropy  (per unit time) of the path measure of the process with respect to the path measure of 
the time-reversed process.  By construction the entropy production rate is nonnegative 
and it vanishes if  and only if the process is reversible.  Crucially, from a numerical point of 
view, the entropy  production rate is an {\em a posteriori} quantity, hence it can be computed
in the course of a simulation as the ergodic average of a certain functional  of the process
(the so-called Gallavotti-Cohen (GC) action functional).  We compute the entropy production
for various numerical schemes such as  explicit Euler-Maruyama and explicit Milstein's for
reversible SDEs with additive or multiplicative noise.  In addition we analyze the entropy
production for the BBK integrator
% and for the exact OU integrator 
for the Langevin equation.  The order (in the time-discretization step $\Delta t$) of the entropy production 
rate provides a tool to classify numerical schemes  in terms of their (discretization-induced) 
irreversibility.  Our results show that the type of the noise critically affects the behavior of the 
entropy production rate. As an example of our results we show that the Euler scheme for {\it multiplicative 
noise} is not an adequate scheme from a reversibility point of view.

\end{abstract}

\begin{resume} Pour un processus de Markov la condition de balance d\'etaill\'ee est \'equivalente 
\`a  la reversibilit\'e du processus par rapport au renversement du temps.  Pour des \'equations 
diff\'erentielles stochastiques,  les sch\'emas de discr\'etisation d\'etruisent en g\'en\'eral cette 
propriet\'e de reversibilit\'e.  En d\'epit d'une vaste litt\'erature  sur l'analyse num\'erique des 
\'equations differentielles stochastiques, leur propriet\'e de stabilit\'e,  les erreurs fortes et/ou faibles, 
les propriet\'es de grandes d\'eviations et \`a long temps, il n'y a pas eu jusqu'\`a maintenant de
r\'esultats quantitatifs sur l'irr\'eversibilit\'e introduite par les approximation num\'eriques. Dans cet article  
nous fournissons de telles estimations, en nous basant sur le taux de production d'entropie,  
inspir\'es par des id\'ees de m\'ecanique statistique hors-\'equilibre. Le taux de production 
d'entropie est, par d\'efinition, l'entropie relative (par unit\'e de temps)  du processus par rapport au 
processus  renvers\'e en temps. Par construction,  le taux de production d'entropie est non-n\'egatif 
et il est z\'ero si et seulement  si le procesus est  r\'eversible.  Crucialement, d'un point de vue num\'erique,
le taux de production d'entropie peut  \^etre evalu\'e  directement comme la moyenne ergodique
d'une certaine fonctionnelle du processus (la fonctionelle de Gallavotti-Cohen), sous des conditions 
d'ergodicit\'e ad\'equates. 
Nous calculons la production d'entropie pour le sch\'ema explicite d'Euler-Maruyama et le 
sch\'ema explicite de Milstein pour des equations diff\'erentielles stochastiques reversibles
avec des bruit additifs ou multiplicatifs.  Nos r\'esultats d\'emontrent que le type de bruit change 
le comportement la production d'entropie de mani\`ere critique. Finalement nous analysons 
la production d'entropie pour le sch\'ema BBK pour l'\'equation de Langevin. 
 \end{resume}



\maketitle

\linenumbers

\section*{Introduction}
\label{intro}
In molecular dynamics algorithms arising in the simulation of systems in materials science, 
chemical engineering,   evolutionary games, computational  statistical mechanics, etc. the steady-
state statistics obtained from numerical  simulations is of great importance \cite{Gardiner:85, vanKampen:06, 
Schlick:02}.   For instance, the free energy of the system or free energy differences as well dynamical 
transitions between metastable states are quantities which are sampled in the stationary regime. 
In addition, physical processes are often modeled at a microscopic level as interactions 
between particles which obey a system of stochastic differential equations (SDE's) 
\cite{Lelievre:10, Gardiner:85}.  To perform steady-state simulations for the sampling of 
desirable observables, the solution of the system of  SDE's  must possess a (unique) ergodic 
invariant measure. The uniqueness of the invariant measure follows from the ellipticity or 
hypoellipticity of the generator of the process together with irreducibility, which 
means that the process can reach at some positive time any open subset of the state space with 
positive probability \cite{Meyn:93,Rey-Bellet:06}.  Under such conditions the distribution 
process converges  to the invariant measure (ergodicity)  which has a smooth density
and the process started in the invariant  measure is stationary, i.e. the distribution 
of the paths of the  processes, is invariant under time-shift.  Many processes of physical origin, 
such as diffusion and adsoprtion/desoprtion of interacting particles, satisfy the condition of 
detailed balance (DB), or equivalently,  time-reversibility, i.e., the distribution of the path of the processes 
are invariant under time-reversal.  It is easy to see that time-reversibility implies stationarity but this a 
strictly stronger condition in general.  The condition of detailed balance often arises from a gradient-like 
behavior of the dynamics or from  Hamiltonian dynamics if the time-reversal includes reversal of the velocities. 
%
%irreducibility \cite{Lelievre:10}. Comprehensively, stationarity means
%that the distribution of a time-shifted path remains the same, reversibility means that
%the probability of observing a path equals the probability of observing the reversed-time path
%while irreducibility means that the process can reach at any positive time any subset of the
%state space with positive probability.
%It is not difficult to show that stationarity and (aperiodic) irreducibility imply the uniqueness of
%the invariant measure\footnote{Invariant measure is also called stationary measure or stationary
%distribution.} as well ergodicity \cite[Thm 17.1.7]{Meyn:93} while irreducibility is a direct
%consequence of the existence of a stationary measure having a positive density with respect
%to Lebesgue measure. Moreover, reversibility implies stationarity while, at equilibrium
%(i.e. stationary) regime, a stochastic process is reversible if and only if Detailed Balance
%(DB) condition is satisfied. 
%%DB is a strong condition that many SDE models with
%%gradient-type drift satisfy.

However, the numerical simulation of  SDE's  necessitates the use of numerical discretization 
schemes. Discretization procedures, except in very special cases,  results in the destruction of 
the DB condition. This affects the approximation process in at least two ways. First,
the invariant measure of the approximation process, if it exists at all, is not known explicitly 
and, second, the time reversibility of the process is lost. Several recent results 
prove the existence of the invariant measure for the discrete-time approximation and 
provide  error estimates  \cite{Bally:96a,Bally:96b,Mattingly:02,Mattingly:10} but, to the best 
of our knowledge, there is no quantitative  assessment of the irreversibility of the 
approximation process.  Of course there exist Metropolized numerical schemes such as MALA
\cite{Roberts:96} and variations thereof which do satisfy the DB condition but they are numerically 
more expensive, especially  in high-dimensional systems, as they require an accept/reject step.
Thus, a quantitative understanding of the lack of reversibility for simpler discretization 
schemes can provide new insights for selecting which schemes are closer to satisfying the DB condition. 

The implications of irreversibility are only partially understood, both from the physical and 
mathematical point of view. These issues have emerged as a main theme in non-equilibrium 
statistical mechanics and it is well-known that irreversibility introduces a stationary current 
(net flow) to the system \cite{Nicolis:77, Schnakenberg:76, Maes:00, Rey-Bellet:11}  but it is unclear how 
this current affects the long-time properties (i.e., the dynamics and large deviations) of 
the process such as exit times, correlation times and phase transitions of metastable states.  
Reversibility is a natural and fundamental property of physical systems and thus, if 
numerical simulation results in the destruction of reversibility, one should carefully {\it quantify
the irreversibility of the approximation process} and we do in this paper {\it using the entropy 
production rate}.  The entropy production rate which is defined as the relative entropy (per unit time)
of the path measure of the process with respect to the path measure of the reversed process is 
widely used in statistical mechanics for the study of non-equilibrium steady states of 
irreversible systems \cite{Gallavotti:95,  Lebowitz:99, Maes:00,Rey-Bellet:11}.  
A fundamental result on the structure of non-equilibrium steady states is the 
Gallavotti-Cohen fluctuation theorem that describes the fluctuations (of large deviations type) 
of the entropy production \cite{Gallavotti:95, Lebowitz:99, Maes:00,Rey-Bellet:11} and this result 
can be viewed as a  generalization of the Kubo-formula and Onsager relations far from equilibrium. 
For our purpose, it is important to note  that the entropy production rate is zero when 
the process is reversible and positive otherwise making {\it entropy production rate a sensible
quantitative  measure of irreversibility}.  Furthermore, if we assume ergodicity of the approximation
process, the entropy production rate equals the time-average of the Gallavotti-Cohen (GC) 
action functional which is defined as the logarithm of the Radon-Nikodym derivative between the 
path measure of the process and the path measure of the reversed process. A key observation
of this paper is that GC action functional is an {\em a posteriori} quantity, hence, it is easily
computable during the simulation making {\it the numerical computation of entropy production
rate tractable}. We show that entropy production is a computable observable that distinguishes
between different numerical schemes in terms of their discretization-induced irreversibility
and as such could allow us to adjust the discretization in the course of the simulation.

We use entropy production to assess the irreversibility of various numerical schemes for 
reversible continuous-time processes. A simple class of reversible processes, yet of great interest, is 
the overdamped Langevin process with gradient-type drift \cite{Gardiner:85, Gillespie:92, 
Lelievre:10}. The discretization of the process is performed using the explicit Euler-Maruyama (EM)
scheme and we distinguish between two different cases depending on the kind of the noise.
In the case of additive noise,  under the assumption of ergodicity of the approximation process 
\cite{Bally:96a,Bally:96b,Mattingly:02,Mattingly:10} we prove that the
entropy production rate is of order $O(\Delta t^2)$ where $\Delta t$ is the time step of the numerical 
scheme. In the case of multiplicative noise,  the results are 
strikingly different. Indeed, under ergodicity assumption, the entropy production
rate for the explicit EM scheme is proved to have a lower positive bound which is independent
of $\Delta t$. Thus irreversibility is {\em not} reduced by adjusting $\Delta t$, as the approximation process 
converges to the continuous-time process. The different behavior of entropy production depending on
the kind of noise is one of the prominent findings of this paper. As a further step in our study, 
we analyze the explicit Milstein's scheme with multiplicative noise (it is the
next higher-order numerical scheme). We prove that the entropy production
rate of Milstein's scheme decreases as time step decreases with order $O(\Delta t)$.

Finally, we compute both analytically and numerically the entropy production rate for a 
discretization scheme for Langevin systems which is another important and widely-used class
of reversible models \cite{Gardiner:85, Lelievre:10}.  The Langevin equation is time-reversible
if addition to reversing time, one reverses the sign of the velocity of all particles.  The noise is 
degenerate but the process is hypo-elliptic and under mild conditions the  Langevin equation 
is ergodic \cite{Talay:02, Mattingly:02, Rey-Bellet:02}.  Our discretization scheme is a
quasi-symplectic splitting scheme also known as BBK integrator  \cite{Brunger:84,Lelievre:10}. 
We rigorously prove, under ergodicity assumption of the  approximation process,  that the entropy rate 
produced by the numerical scheme for the  Langevin process with additive  noise is of order $O(\Delta t)$, 
hence, in terms of irreversibility it  is an acceptable integration scheme.

The paper is organized in four sections. In Section~\ref{entropy:prod:section} we recall some 
basic facts about  reversible processes and define the entropy production. Moreover 
we give the basic assumptions necessary for our proofs, namely, the ergodicity of both 
continuous-time and discrete-time approximation process.  In Section~
\ref{overdamped:Langevin:sec} we compute the entropy production rate for reversible 
overdamped Langevin processes. The section is split into three subsections for the additive and 
multiplicative noise for the Euler and Milstein schemes. In Section~\ref{Langevin:sec} we compute the entropy production rate for the 
reversible (up to momenta flip) Langevin process using the BBK integrator. Conclusions and future 
extensions of the current work are summarized in the fourth and final Section.

\section{Reversibility, Gallavotti-Cohen Action Functional, and Entropy Production}
\label{entropy:prod:section}
Let us consider a $d$-dimensional system of SDE's written as
\begin{equation}
dX_t = a(X_t)dt + b(X_t) dB_t
\label{general:sde}
\end{equation}
where $X_t\in\mathbb R^d$ is a diffusion Markov process, %driven by (\ref{general:sde}),
$a:\mathbb R^d\rightarrow\mathbb R^d$ is the drift vector,
$b:\mathbb R^d\rightarrow\mathbb R^{d\times m}$ is the diffusion matrix, 
and $B_t\in\mathbb R^m$ is a standard $m$-dimensional
Brownian motion.  We will always assume that $a$ and $b$ are sufficiently smooth and satisfy 
suitable  growth conditions and/or dissipativity conditions at infinity  to ensure the existence of 
global solutions.  The generator of the diffusion process is defined by 
\begin{equation}
\mathcal Lf = \sum_{i=1}^d a_i\frac{\partial f}{\partial x_i} + \frac{1}{2}\sum_{i,j=1}^d (bb^T)_{i,j}\frac{\partial^2 f}{\partial x_i\partial x_j}.
\label{sde:generator}
\end{equation}
for smooth test functions $f$. 
% which are twice continuously differentiable and with bounded 
%derivatives up to second order.
We assume that the process $X_t$ has a (unique) invariant measure $\mu(dx)$, 
and that it satisfies the Detailed Balance (DB) condition, i.e., its generator is
symmetric in the Hilbert space $L^2(\mu)$: 
\begin{equation}
<\mathcal Lf,g>_{L^2(\mu)} = <f,\mathcal Lg>_{L^2(\mu)}
\label{generator:symmetry:cond}
\end{equation}
for suitable smooth test functions $f,g$. 

A Markov process $X_t$ is said to be time-reversible if for any $n$ and sequence of times 
$t_1< \cdots<t_n$ the finite dimensional distributions of  $(X_{t_1},...,X_{t_n})$
and of $(X_{t_n},...,X_{t_1})$ are identical. More formally, let ${\bf P}^\rho_{[0,t]}$ denote
the path measure of the process $X_t$ on the time-interval $[0,t]$  with $X_0 \sim \rho$. 
Let $\Theta$ denote the time reversal, i.e. $\Theta$ acts on a path $\{X_s\}_{0\le s \le t}$ has 
\begin{equation}
(\Theta X)_s  \,=\, X_{t-s}
\end{equation}
Then reversibility is equivalent to  ${\bf P}^\mu_{[0,t]} \,=\, {\bf P}^\mu_{[0,t]} \circ \Theta$
and  it is well-known that a stationary\footnote{Stationarity is equivalent to starting
the process $X_t$ from its invariant measure, i.e., $X_0\sim\mu$.} process
which satisfies the DB condition is time-reversible.

The condition of reversibility can be also expressed in terms of relative entropy as follows. 
Recall that for two probability measure $\pi_1, \pi_2$ on some measurable space, the relative 
entropy of $\pi_1$ with respect to $\pi_2$ is given by  $R(\pi_1|\pi_2) \equiv \int d\pi_1 \log\frac{d
\pi_1}{d\pi_2}$ if $\pi_1$ is absolutely continuous with respect to $\pi_2$ and $+\infty$ otherwise. 
The relative entropy is nonnegative, $R(\pi_1|\pi_2)\ge 0$ and  $R(\pi_1|\pi_2)= 0$ if and only if 
$\pi_1=\pi_2$. The entropy production rate of a  Markov process $X_t$ is defined by 
\begin{equation}
EP_{cont} := \lim_{t\rightarrow\infty} \frac{1}{t} R({\bf P}^\rho_{[0,t]}|{\bf P}^\rho_{[0,t]} \circ \Theta)
= \lim_{t\rightarrow\infty} \frac{1}{t}\int d{\bf P}^\rho_{[0,t]} \log\frac{d{\bf P}^\rho_{[0,t]}}{d{\bf P}^\rho_{[0,t]}\circ \Theta}
\label{entropy:prod:definition}
\end{equation}
If $X_t$ satisfies DB and $X_0 \sim \mu$ then $R({\bf P}^\mu_{[0,t]}|{\bf P}^\mu_{[0,t]} \circ \Theta)$
is identically $0$ for all $t$ and the entropy production rate is $0$. Note that if $X_0 \sim \rho 
\not=\mu$ then $R({\bf P}^\rho_{[0,t]}|{\bf P}^\rho_{[0,t]} \circ \Theta)$ is a boundary term, in the sense that it is  
$O(1)$ and so the entropy rate vanishes in this case in the large time limit  (under suitable ergodicity 
assumptions).  Conversely when $EP_{cont}\neq 0$ the process is  truly irreversible.  The entropy 
production rate  for  Markov processes and stochastic differential  equations  is discussed in more detail in  \cite{Lebowitz:99,Maes:00}.  
 
%
%
%Another way to prove the reversibility of a stationary process is to show
%that the entropy production rate of the process is zero. Indeed, denoting by $P_{0t}$
%the path measure of the stationary process in the time window $[0,t]$ and by $P_{t0}$
%the path measure of the reversed-time process, entropy production rate is defined
%as the time-mean relative entropy between the two processes
%\begin{equation}
%EP_{cont} := \lim_{t\rightarrow\infty} \frac{1}{t} R(P_{0t}|P_{t0})
%= \lim_{t\rightarrow\infty} \frac{1}{t} \mathbb E_P[\log\frac{dP_{0t}}{dP_{t0}}]
%= \lim_{t\rightarrow\infty} \frac{1}{t}\int dP \log\frac{dP_{0t}}{dP_{t0}}
%\label{entropy:prod:definition}
%\end{equation}
%where $\frac{dP_{0t}}{dP_{t0}}$ is the Radon-Nikodym derivative between the two
%processes while $P$ is the path space measure of the process having started from its
%invariant measure. Obviously, when the process is time reversible, the path measures $P_{0t}$
%and $P_{t0}$ are equal, thus, the Radon-Nikodym derivative equals to 1 $P-a.s.$ and,
%consequently, $EP_{cont}=0$. Conversely, when $EP_{cont}\neq 0$ the process is irreversible.

Let us consider a numerical integration scheme for the SDE (\ref{general:sde}) which has
the general form
\begin{equation}
x_{i+1} = F(x_i,\Delta t, \Delta W_i) \ \ \ i = 1,2,...
\label{sde:num:scheme}
\end{equation}
Here $x_i \in\mathbb R^d$ is a discrete-time continuous state-space Markov process,
$\Delta t$ is the time-step and  $\Delta W_i \in\mathbb R^m,\ i=1,2,...$
are i.i.d. Gaussian  random variables with mean $0$ and variance $\Delta t I_m$. 
We will assume that the Markov process $x_i$ has transition probabilities which are
absolutely continuous with respect to Lebesgue measure with everywhere positive 
densities  $\Pi(x_i,x_{i+1}) := \Pi_{F(x,\Delta t, \Delta W)}(x_{i+1}|x_i)$
and we also assume that  $x_i$ has a invariant measure which we denote $\bar{\mu}(dx)$ 
and which is then unique and  has a density with respect to Lebesgue.  In general the 
invariant measure for $X_t$ and $x_i$ differ,  $\mu \not= \bar{\mu}$ and $x_i$ does not satisfy a 
DB condition.  Note also that the very existence of $\bar{\mu}$ is 
not guaranteed in general.  Results on the existence of $\bar{\mu}$ do exist however and 
typically require  that the SDE is elliptic or hypoellitptic and that the state space of $X_t$ is 
compact or that a  global Lipschitz condition on the drift  holds 
\cite{Bally:96a,Bally:96b,Mattingly:02,Mattingly:10}. 

Proceeding as in the continuous case we introduce an entropy production rate for the Markov 
process $x_i$. Let us assume that the process starts from some distribution $\rho(x)dx$, then the
finite dimensional distribution on the time window $[0,t]$ where $t=n\Delta t$ 
is given by 
\begin{equation}
\bar{{\bf P}}_{[0,t]}(dx_0,\cdots,  dx_n) = \rho(x_0)\Pi(x_0,x_1)\cdots\Pi(x_{n-1},x_n) dx_0 
\cdots dx_n \,.
\end{equation}
For the time reversed path $\Theta (x_0, \cdots x_n) \,=\, (x_n, \cdots, x_0)$ we have then
\begin{equation}
\bar{{\bf P}}_{[0,t]}\circ \Theta(dx_0,...,dx_n) = \rho(x_n)\Pi(x_n,x_{n-1})\cdots\Pi(x_1,x_0) dx_0 
\cdots dx_n 
\end{equation}
and the Radon-Nikodym derivative takes the form 
\begin{equation}
\frac{d \bar{{\bf P}}_{[0,t]}}{d \bar{\bf P}_{[0,t]}\circ \Theta} = \exp(W(t)) \frac{\rho(x_0)}{\rho(x_n)}
\label{Radon:Nikodym:der}
\end{equation}
where $W(t)$ is the Gallavotti-Cohen (GC) action functional
given by
\begin{equation}
W(t) = W(n;\Delta t) := \sum_{i=0}^{n-1} \log\frac{\Pi(x_i,x_{i+1})}{\Pi(x_{i+1},x_i)} \,.
\label{GC:action:func}
\end{equation}
Note that $W(t)$ is an additive functional of the paths and thus  if $x_i$ is ergodic,  
by the  ergodic theorem the following limit exists 
\begin{equation}
EP(\Delta t) = \lim_{t\rightarrow\infty} \frac{1}{t} W(t) = \lim_{n\rightarrow\infty} \frac{1}{n\Delta t} W(n;\Delta t) \ \ \ \bar{P}-a.s..
\label{entr:prod:GCAF}
\end{equation}
We call the quantity $EP(\Delta t)$  the entropy production rate associated to the 
numerical scheme. Note that we have, almost surely,  
\begin{equation}
EP(\Delta t) \,=\, \frac{1}{\Delta t} \lim_{n \to \infty} \frac{1}{n} 
\sum_{i=0}^{n-1} \log\frac{\Pi(x_i,x_{i+1})}{\Pi(x_{i+1},x_i)} \,=\, 
\frac{1}{\Delta t}
\int\int  \bar{\mu}(x) \Pi(x,y) \log\frac{\Pi(x,y)}{\Pi(y,x)} \,dx dy   
\end{equation}
and for concrete numerical schemes we will compute fairly explicitly the entropy production
in the next sections.  Since we are 
interested in the ergodic average we will systematically omit 
boundary  terms  which do not contribute  to ergodic averages and we will use the notation 
\begin{equation}\label{equivalent}
W_1(t)\dot{=}W_2(t) \quad {\rm ~if~} \quad   \lim_{t \to \infty } \frac{1}{t} 
(W_1(t) -W_2(t))=0\,.
\end{equation} 
For example we have 
\begin{equation}
W(t) \,\dot{=}\, \log \frac{d \bar{{\bf P}}_{[0,t]}}{d \bar{\bf P}_{[0,t]}\circ \Theta} \,.
\end{equation} 
Note also that using (\ref{entr:prod:GCAF}) and (\ref{GC:action:func}),
entropy production rate is tractable numerically and it can be easily calculated
``on-the-fly" once the transition probability density function $\Pi(\cdot,\cdot)$ is provided.

In  the following sections we investigate the behavior of the entropy production rate for different 
discretization schemes of various reversible processes in the stationary regime.
However, before proceeding with our analysis, let us state formally the basic assumptions 
necessary for our results to apply.
\begin{ssmptn} \label{basic:assumption} We have 
\begin{itemize}
\item The drift $a$ and the diffusion $b$ in \eqref{general:sde} as well as the vector $F$ in 
\eqref{sde:num:scheme} are $C^\infty$ and all their derivatives 
have at most polynomial growth at infinity.   
\item  The generator $\mathcal L$ is elliptic or hypo-elliptic, in particular the transition 
probabilities and the invariant measure (if it exists) are absolutely continuous with respect to
Lebesgue with smooth densities. We assume that $x_t$ is ergodic, i.e. every open set can 
be reached with positive probability starting from any point.  For the discretized scheme we assume 
that $x_i$ has smooth  everywhere positive transition probabilities.   
\item Both the continuous-time process $X_t$ and discrete-time process $x_i$ are ergodic with
unique invariant measures $\mu$ and $\bar{\mu}$, respectively.  Furthermore for sufficiently small 
$\Delta t$ we have 
\begin{equation}
\left| \mathbb E_\mu[f] - \mathbb E_{\bar{\mu}}[f]\right| = O(\Delta t)
\label{stat:average:estimate}
\end{equation}
for functions $f$ which are $C^\infty$ with at most polynomial growth at infinity. 
\end{itemize}
\end{ssmptn}

Notice that inequality (\ref{stat:average:estimate}) is an error estimate for the invariant measures
of  the processes $X_t$ and $x_i$. The rate of convergence in terms of
$\Delta t$ depends on the particular numerical scheme \cite{Talay:90a, Mattingly:10}. Ergodicity
results for (numerical) SDEs can be found in \cite{Roberts:96, Mattingly:02, Talay:02,
Talay:90a, Talay:90b, Bally:96a, Bally:96b, Mattingly:10, Khasminskii:10}. 
For instance, if both drift term $a(x)$ and diffusion term $b(x)$ have bounded derivatives of 
any order, the covariance matrix $(bb^T)(x)$ is elliptic for all $x\in\mathbb R^d$ and there is a 
compact set outside of which holds $x^Ta(x)<-C|x|^2$ for all $x\in\mathbb R^d$ (Lyapunov 
exponent) then it was shown in \cite{Talay:90a} that the continuous-time process as well both Euler 
and Milstein numerical schemes are ergodic and error estimate (\ref{stat:average:estimate}) holds.
Another less restrictive example where ergodicity properties were proved is for SDE systems
with degenerate noise and particularly for Langevin processes \cite{Mattingly:02, Talay:02}.
Again, a Lyapunov functional is the key assumption in order to handle the stochastic process
at the infinity. More recently, Mattingly et al. \cite{Mattingly:10} showed ergodicity for
SDE-driven processes restricted on a torus as well their discretizations utilizing only
the assumptions of ellipticity or hypoellipticity and the assumption of local 
Lipschitz continuity for both drift and diffusion terms.

\section{Entropy Production for  Overdamped Langevin Processes}
\label{overdamped:Langevin:sec}
The overdamped Langevin process, $X_t\in\mathbb R^d$, is the solution of the following
system of SDE's
\begin{equation}
\begin{aligned}
&dX_t = - \frac{1}{2}\Sigma(X_t)\nabla V(X_t)dt + \frac{1}{2}\nabla \Sigma(X_t)dt + \sigma(X_t) dB_t 
%&X_0 \sim \text{invariant measure}
\label{multi:sde:cont}
\end{aligned}
\end{equation}
where $V:\mathbb R^d \rightarrow \mathbb R$ is a smooth potential function,
$\sigma:R^d \rightarrow \mathbb R^{d\times m}$ is the diffusion matrix,
$\Sigma:=\sigma\sigma^T:R^d \rightarrow \mathbb R^{d\times d}$ is the covariance matrix
and $B_t$ is a standard 
$m$-dimensional Brownian motion.  We assume from now on that $\Sigma(x)$ is invertible
for any $x$ so that the process is elliptic.  It is straightforward to show that the generator of the 
process $X_t$ satisfies the DB condition (\ref{generator:symmetry:cond}) with invariant measure
\begin{equation}
\mu(dx) = \frac{1}{Z} \exp(-V(x)) dx
\label{inv:measure:overdamed}
\end{equation}
where $Z=\int_{\mathbb R^d} \exp(-V(x)) dx$ is the normalization constant and thus if $X_0 \sim \mu$
then the Markov process $X_t$ is reversible. 
%
%Moreover, the
%satisfaction of the DB condition plus the fact that the process was started from the
%invariant measure (i.e., equilibrium regime) provides the time reversibility property of the process.

The explicit Euler-Maruyama (EM) scheme for numerical integration of
(\ref{multi:sde:cont}) is given by
\begin{equation}
x_{i+1} = x_i -\frac{1}{2}\Sigma(x_i)\nabla V(x_i)\Delta t + \frac{1}{2}\nabla \Sigma(x_i)\Delta t + \sigma(x_i) \Delta W_i 
\label{multi:sde:explEM}
\end{equation}
with $\Delta W_i \sim N(0, \Delta t I_m)$, $i=1,2,...$ are $m$-dimensional iid Gaussian random
variables. The process $x_i$ is a discrete-time Markov process with transition probability 
density given by
\begin{equation}
\begin{aligned}
\Pi(x_i,x_{i+1}) = \frac{1}{Z(x_i)} &\exp\left(\frac{1}{2\Delta t} (\Delta x_i + \frac{1}{2}\Sigma(x_i)\nabla V(x_i)\Delta t-\frac{1}{2}\nabla \Sigma(x_i)\Delta t)^T \right.\\
&\left.\Sigma^{-1}(x_i) (\Delta x_i + \frac{1}{2}\Sigma(x_i)\nabla V(x_i)\Delta t-\frac{1}{2}\nabla \Sigma(x_i)\Delta t) \right)
\label{trans:prob:explEM:multi:sde}
\end{aligned}
\end{equation}
where $\Delta x_i = x_{i+1}-x_i$ and $Z(x_i)=(2\pi)^{m/2}|\det \Sigma(x_i)|^{1/2}$ is the
normalization constant for the multidimensional Gaussian distribution. The 
following lemma provides the GC action functional for the explicit EM time-discretization
scheme of the overdamped Langevin process.

\begin{lmm}
Assume that $\det\Sigma(x)\neq0\ \forall x\in\mathbb R^d$.  Then the GC action functional of 
the process $x_i$ solving (\ref{multi:sde:explEM}) is
\begin{equation}
\begin{aligned}
W(n;\Delta t) &\dot{=} - \frac{1}{2}\sum_{i=0}^{n-1} \Delta x_{i}^T[\nabla V(x_{i+1})+\nabla V(x_i)]
+ \frac{1}{2}  \sum_{i=0}^{n-1} \Delta x_{i}^T[\Sigma^{-1}(x_{i+1})\nabla \Sigma(x_{i+1}) + \Sigma^{-1}(x_i)\nabla \Sigma(x_i)] \\
&+\frac{1}{2\Delta t} \sum_{i=0}^{n-1} \Delta x_i^T\left[\Sigma^{-1}(x_{i+1})-\Sigma^{-1}(x_i)\right]\Delta x_i 
\label{action:func:multi}
\end{aligned}
\end{equation}
where $\dot{=}$ means equality up to boundary terms, as defined in \eqref{equivalent}.
\label{overdamped:Langevin:lemma}
\end{lmm}

\begin{proof}
The assumption for non-zero determinant
is imposed so that the transition probabilities and hence the GC action functional
are non-singular.
The proof is then a straightforward computation using (\ref{trans:prob:explEM:multi:sde}) and
(\ref{GC:action:func}). 
\begin{equation*}
\begin{aligned}
&W(n;\Delta t) := \sum_{i=0}^{n-1} \left[\log\Pi(x_i,x_{i+1}) - \log\Pi(x_{i+1},x_i) \right]
= \sum_{i=0}^{n-1} \left[\log Z(x_{i+1}) - \log Z(x_i) \right]  \\
&-\frac{1}{2\Delta t} \sum_{i=0}^{n-1} \left[ (\Delta x_i + \frac{1}{2}\Sigma(x_i)\nabla V(x_i)\Delta t-\frac{1}{2}\nabla \Sigma(x_i)\Delta t)^T
\Sigma^{-1}(x_i) (\Delta x_i + \frac{1}{2}\Sigma(x_i)\nabla V(x_i)\Delta t-\frac{1}{2}\nabla \Sigma(x_i)\Delta t) \right. \\
&\left.- (-\Delta x_i + \frac{1}{2}\Sigma(x_{i+1})\nabla V(x_{i+1})\Delta t-\frac{1}{2}\nabla \Sigma(x_{i+1})\Delta t)^T
\Sigma^{-1}(x_{i+1}) (-\Delta x_i + \frac{1}{2}\Sigma(x_{i+1})\nabla V(x_{i+1})\Delta t-\frac{1}{2}\nabla \Sigma(x_{i+1})\Delta t) \right] \\
&\dot{=}-\frac{1}{2\Delta t} \sum_{i=0}^{n-1} \left[ \Delta x_i^T\Sigma^{-1}(x_i)\Delta x_i\right.
+ \frac{1}{4}\nabla V(x_i)^T\Sigma(x_i)\nabla V(x_i)\Delta t^2 + \frac{1}{4}\nabla \Sigma(x_i)^T \Sigma^{-1}(x_i)\nabla \Sigma(x_i)\Delta t^2 \\
&+ \Delta x_i^T\nabla V(x_i)\Delta t - \Delta x_i^T\Sigma^{-1}(x_i)\nabla \Sigma(x_i)\Delta t
-\frac{1}{2}\nabla V(x_i)^T\nabla \Sigma(x_i)\Delta t^2 \\
&-\Delta x_i^T\Sigma^{-1}(x_{i+1})\Delta x_i
- \frac{1}{4}\nabla V(x_{i+1})^T\Sigma(x_{i+1})\nabla V(x_{i+1})\Delta t^2 - \frac{1}{4}\nabla \Sigma(x_{i+1})^T \Sigma^{-1}(x_{i+1})\nabla \Sigma(x_{i+1})\Delta t^2 \\
&\left.+ \Delta x_i^T\nabla V(x_{i+1})\Delta t - \Delta x_i^T\Sigma^{-1}(x_{i+1})\nabla \Sigma(x_{i+1})\Delta t
+\frac{1}{2}\nabla V(x_{i+1})^T\nabla \Sigma(x_{i+1})\Delta t^2\right] \\
&\dot{=} -\frac{1}{2\Delta t} \sum_{i=0}^{n-1} \Delta x_i^T\left[\Sigma^{-1}(x_i)-\Sigma^{-1}(x_{i+1})\right]\Delta x_i
- \frac{1}{2}\sum_{i=0}^{n-1} \Delta x_{i}^T[\nabla V(x_{i+1})+\nabla V(x_i)] \\
&+ \frac{1}{2}  \sum_{i=0}^{n-1} \Delta x_{i}^T[\Sigma^{-1}(x_{i+1})\nabla \Sigma(x_{i+1}) + \Sigma^{-1}(x_i)\nabla \Sigma(x_i)]
\end{aligned}
\end{equation*}
where all the terms of the general form $G(x_i)-G(x_{i+1})$ in the sums were cancelled out
since they form telescopic sums which become boundary terms.
\end{proof}
Three important remarks can readily be made from the above computation. 

\begin{rmrk}
The numerical computation of entropy production rate as the time-average of the 
GC action functional on the path space (i.e., based on 
%(\ref{entropy:prod:discrete:process})
\eqref{Radon:Nikodym:der}) at first sight seems  computationally intractable due to the large 
dimension of the path space. However,  due to ergodicity, the numerical computation of the 
entropy production can be 
performed as a time-average based on (\ref{entr:prod:GCAF}) and (\ref{action:func:multi}) for large 
$n$. Additionally, this computation can be done for free and ``on-the-fly" since the quantities involved are 
already computed in the simulation of the process.  The numerical entropy production rate
shown in the following figures is computed using this approach.
% Moreover, someone wants to create an adaptive scheme that is reversible
%or almost reversible, a starting point could be the further study and understanding
%of the GC action functional.
\end{rmrk}

\begin{rmrk}
\label{over:Langevin:remark}
It was shown in \cite{Maes:00} that the GC action functional of the {\it continuous-time} process
driven by (\ref{multi:sde:cont}) equals the Stratonovich integral
\begin{equation}
W_{cont}(t) = -\int_0^t \nabla V(X_s) \circ dX_s = V(x_0) - V(x_t) 
\end{equation}
which reduces to a boundary term as expected.  This functional has the discretization 
\begin{equation}
W_{cont}(t)  \approx \frac{1}{2} \sum_{i=0}^{n-1} \Delta x_{i}^T[\nabla V(x_{i+1})+\nabla V(x_i)]
\end{equation}
and this is exactly the first term in the GC action functional $W(n;\Delta t)$ 
for the explicit EM approximation  process (see (\ref{action:func:multi})). 
%we observe that the discretized $W_{cont}(t)$ is contained in the GC action functional
%of the discrete-time approximation process $W(n;\Delta t)$. 
However, the discretization scheme introduces two additional terms to the GC action 
functional which may  greatly affect the asymptotic behavior of entropy production as 
$\Delta t$ goes to zero, as we demonstrate in Section~\ref{multiplicative}. Notice that when the noise is additive, i.e., when the diffusion matrix is 
constant, then these two  additional terms vanish and  taking the limit $\Delta t\rightarrow0$,
the GC action functional $W(n;\Delta t)$, if exists, becomes the Stratonovich
integral $W_{cont}(t)$ which is a boundary term.  
%Hence, the entropy production
%for the additive noise case is zero and the continuous-time process $X_t$ is
%reversible. This is an alternative approach to show that the process $X_t$ is reversible
%without proving the satisfaction of detailed balance condition which necessitates the
%knowledge of the invariant measure.
\end{rmrk}

\begin{rmrk}
The GC action functional $W(n;\Delta t)$ consists of three terms (see (\ref{action:func:multi})), each of 
which  stems from a particular term in the SDE. Thus, each term in the SDE
contributes to the entropy production functional  a component  which is totally decoupled
to the other terms. The reason for  this decomposition lies  in the particular form of the
transition probabilities for the explicit EM scheme which are exponentials with
quadratic argument. This feature can be exploited for the study of entropy
production of numerical schemes for processes with irreversible dynamics.
Indeed, if a non-gradient term of the form $a(X_t)dt$ is added to the drift of (\ref{multi:sde:cont}), the process is 
irreversible and its GC action functional is not anymore a boundary term  and is given by 
\cite{Maes:00}
\begin{equation}
W_{cont}(t) \dot{=} -\int_0^t \Sigma^{-1}(X_t) a(X_t) \circ dX_t \approx
\frac{1}{2} \sum_{i=0}^{n-1} \Delta x_{i}^T[\Sigma^{-1}(x_i) a(x_i) + \Sigma^{-1}(x_{i+1}) a(x_{i+1})]
\end{equation}
On the other hand,
due to the separation property of the explicit EM scheme, the GC action functional
of the discrete-time approximation process $W(n;\Delta t)$ has the additional term
\begin{equation}
\frac{1}{2} \sum_{i=0}^{n-1} \Delta x_{i}^T[\Sigma^{-1}(x_i) a(x_i) + \Sigma^{-1}(x_{i+1}) a(x_{i+1})].
\end{equation}
Evidently, the discretization of $W_{cont}(t)$ equals the additional term of the GC functional 
$W(n;\Delta t)$. Thus, GC action functional $W(n;\Delta t)$ is decomposed  into two components, one 
stemming from the irreversibility  of the continuous-time process and another one stemming from the 
irreversibility of the discretization procedure.
\end{rmrk}

\subsection{Entropy Production for the Additive Noise Case}
An important special case of (\ref{multi:sde:cont}) is the case of additive 
noise, i.e., when the covariance matrix does not depend in the
process, $\Sigma(x) \equiv \Sigma$. In this case, the SDE system becomes
\begin{equation}
\begin{aligned}
&dX_t = -\frac{1}{2}\Sigma\nabla V(X_t)dt + \sigma dB_t \\
&X_0 \sim \mu
\label{add:sde:cont}
\end{aligned}
\end{equation}
and the GC action functional is simply given by 
\begin{equation}
W(n;\Delta t) \dot{=} -\frac{1}{2}\sum_{i=0}^{n-1} \Delta x_{i}^T[\nabla V(x_{i+1})+\nabla V(x_i)]
\label{action:func:add}
\end{equation}

In this section we prove an upper bound for the entropy production of the explicit EM scheme. 
The proof uses several lemmas stated and proved in Appendix~\ref{app:a}.

\begin{thrm}\label{add:noise:theorem}
Let Assumption~\ref{basic:assumption} hold. Assume also that the potential
function $V$ has bounded fifth-order derivative and that the covariance matrix $\Sigma$ is invertible. 
Then, for sufficiently small $\Delta t$, there exists $C=C(V,\Sigma)>0$ such that
\begin{equation}
EP(\Delta t) \leq C\Delta t^2
\label{entr:prod:add}
\end{equation}
%\label{add:noise:theorem}
\end{thrm}

\begin{proof}
%Mattingly et. al \cite{Mattingly:10} (Theorem 5.3) showed for the explicit Euler scheme
%on a torus that if covariance matrix is positive definite then Assumption~\ref{basic:assumption}(ii)
%is satisfied for all test functions in $W^{2,\infty}$. Moreover, they stated that by suitable
%handling of the drift term (globally Lipschitz) as in \cite{Mattingly:02}, they should be
%able to extend the theorem to $\mathbb R^d$.
Utilizing the generalized trapezoidal rule (\ref{gen:trapez:rule}) for $k=3$,
the GC action function is rewritten as
\begin{equation}
\begin{aligned}
W(n;\Delta t) &\dot{=} - \frac{1}{2} \sum_{i=0}^{n-1} \Delta x_{i}^T[\nabla V(x_{i+1})+\nabla V(x_i)] \\
&= \sum_{i=0}^{n-1} \left\{-(V(x_{i+1})-V(x_i)) + \sum_{|\alpha|=3} C_\alpha [D^\alpha V(x_{i+1})+D^\alpha V(x_i)]\Delta x_i^\alpha \right. \\
&\left.+ \sum_{|\alpha|=1,3,5} \sum_{|\beta|=5-|\alpha|} B_\beta [R_\alpha^\beta(x_i,x_{i+1}) + R_\alpha^\beta(x_{i+1},x_i)]\Delta x_i^{\alpha+\beta} \right\} \\
&\dot{=} \sum_{i=0}^{n-1} \sum_{|\alpha|=3} C_\alpha [D^\alpha V(x_{i+1})+D^\alpha V(x_i)]\Delta x_i^\alpha \\
&+ \sum_{i=0}^{n-1} \sum_{|\alpha|=1,3,5} \sum_{|\beta|=5-|\alpha|} B_\beta [R_\alpha^\beta(x_i,x_{i+1}) + R_\alpha^\beta(x_{i+1},x_i)]\Delta x_i^{\alpha+\beta}\, .
\end{aligned}
\end{equation}
Applying, once again, Taylor series expansion to $D^\alpha V(x_{i+1})$, the GC action functional
becomes
\begin{equation}
\begin{aligned}
W(n;\Delta t) &\dot{=} \sum_{i=0}^{n-1} \left\{\sum_{|\alpha|=3} 2C_\alpha D^\alpha V(x_i) \Delta x_i^\alpha
+ \sum_{|\alpha|=3} C_\alpha \sum_{|\beta|=1} D^{\alpha+\beta} V(x_i)\Delta x_i^{\alpha+\beta} \right\} \\
&+ \sum_{i=0}^{n-1} \sum_{|\alpha|=1,3,5} \sum_{|\beta|=5-|\alpha|} \bar{R}_\alpha^\beta(x_i,x_{i+1}) \Delta x_i^{\alpha+\beta}
\end{aligned}
\end{equation}
where $\bar{R}_\alpha^\beta(x_i,x_{i+1}) = B_\beta [R_\alpha^\beta(x_i,x_{i+1}) + R_\alpha^\beta(x_{i+1},x_i)]
+\mathbbm 1_{|\alpha|=3} R_\beta^\alpha(x_i,x_{i+1})$.
Moreover, expanding $\Delta x_i^\alpha$ using the multi-binomial formula
\begin{equation}
\Delta x_i^\alpha = (-\frac{1}{2}\Sigma\nabla V(x_i)\Delta t + \sigma\Delta W_i)^\alpha
= \sum_{\nu\leq\alpha} \binom{\alpha}{\nu} (-\frac{1}{2}\Sigma\nabla V(x_i)\Delta t)^\nu (\sigma\Delta W_i)^{\alpha-\nu}\, .
\end{equation}
Then, the GC action functional becomes
\begin{equation}
\begin{aligned}
&W(n;\Delta t) \dot{=} 2\sum_{i=0}^{n-1} \sum_{|\alpha|=3} \sum_{\nu\leq\alpha}
C_\alpha \binom{\alpha}{\nu} D^\alpha V(x_i) (-\frac{1}{2}\Sigma\nabla V(x_i)\Delta t)^\nu (\sigma\Delta W_i)^{\alpha-\nu} \\
&+ \sum_{i=0}^{n-1}\sum_{|\alpha|=3} \sum_{|\beta|=1} \sum_{\nu\leq\alpha+\beta}
C_\alpha \binom{\alpha+\beta}{\nu} D^{\alpha+\beta} V(x_i) (-\frac{1}{2}\Sigma\nabla V(x_i)\Delta t)^\nu (\sigma\Delta W_i)^{\alpha+\beta-\nu} \\
&+ \sum_{i=0}^{n-1} \sum_{|\alpha|=1,3,5} \sum_{|\beta|=5-|\alpha|}\sum_{\nu\leq\alpha+\beta}
\binom{\alpha+\beta}{\nu} \bar{R}_\alpha^\beta(x_i,x_{i+1})
(-\frac{1}{2}\Sigma\nabla V(x_i)\Delta t)^\nu (\sigma\Delta W_i)^{\alpha+\beta-\nu}\, .
\end{aligned}
\end{equation}
From (\ref{entr:prod:GCAF}), the entropy production rate is the time-averaged
GC action functional as $n\rightarrow\infty$. Thus,
\begin{equation}
\begin{aligned}
&EP(\Delta t) = \lim_{n\rightarrow\infty} \frac{W(n;\Delta t)}{n\Delta t} \\
&= \frac{2}{\Delta t} \sum_{|\alpha|=3} \sum_{\nu\leq\alpha} C_\alpha \binom{\alpha}{\nu}
\lim_{n\rightarrow\infty} \frac{1}{n} \sum_{i=0}^{n-1} D^\alpha V(x_i) (-\frac{1}{2}\Sigma\nabla V(x_i)\Delta t)^\nu (\sigma\Delta W_i)^{\alpha-\nu} \\
&+ \frac{1}{\Delta t}\sum_{|\alpha|=3} \sum_{|\beta|=1} \sum_{\nu\leq\alpha+\beta} C_\alpha \binom{\alpha+\beta}{\nu} 
\lim_{n\rightarrow\infty} \frac{1}{n} \sum_{i=0}^{n-1} D^{\alpha+\beta} V(x_i) (-\frac{1}{2}\Sigma\nabla V(x_i)\Delta t)^\nu (\sigma\Delta W_i)^{\alpha+\beta-\nu} \\
&+ \frac{1}{\Delta t} \sum_{|\alpha|=1,3,5} \sum_{|\beta|=5-|\alpha|}\sum_{\nu\leq\alpha+\beta} \binom{\alpha+\beta}{\nu}
\lim_{n\rightarrow\infty} \frac{1}{n} \sum_{i=0}^{n-1} \bar{R}_\alpha^\beta(x_i,x_{i+1})
(-\frac{1}{2}\Sigma\nabla V(x_i)\Delta t)^\nu (\sigma\Delta W_i)^{\alpha+\beta-\nu}\, .
\label{entropy:prod:add1}
\end{aligned}
\end{equation}
The ergodicity of $x_i$ as well the Gaussianity of $\Delta W_i$ guarantees that the first two limits
in the entropy production formula exist. Additionally, the residual terms, $\bar{R}_\alpha^\beta(x_i,x_{i+1})$, are
bounded due to the assumption on bounded fifth-order derivative of $V$, hence, the third limit also exists. Note here
that this assumption could be changed by assuming boundedness of a higher order derivative and performing
a higher-order Taylor expansion. Appendix~\ref{app:a} gives rigorous proofs of these ergodicity statements.
Hence, 
\begin{equation}
\begin{aligned}
&EP(\Delta t) = \frac{2}{\Delta t} \sum_{|\alpha|=3} \sum_{\nu\leq\alpha} C_\alpha \binom{\alpha}{\nu}
\mathbb E_{\bar{\mu}} [D^\alpha V(x) (-\frac{1}{2}\Sigma\nabla V(x)\Delta t)^\nu] \mathbb E_{\rho} [(\sigma y)^{\alpha-\nu}] \\
&+ \frac{1}{\Delta t}\sum_{|\alpha|=3} \sum_{|\beta|=1} \sum_{\nu\leq\alpha+\beta} C_\alpha \binom{\alpha+\beta}{\nu} 
\mathbb E_{\bar{\mu}} [D^{\alpha+\beta} V(x) (-\frac{1}{2}\Sigma\nabla V(x)\Delta t)^\nu] \mathbb E_{\rho} [(\sigma y)^{\alpha+\beta-\nu}] \\
&+ \frac{1}{\Delta t} \sum_{|\alpha|=1,3,5} \sum_{|\beta|=5-|\alpha|}\sum_{\nu\leq\alpha+\beta} \binom{\alpha+\beta}{\nu}
\mathbb E_{\bar{\mu}\times\rho} [\bar{R}_\alpha^\beta(x,y) (-\frac{1}{2}\Sigma\nabla V(x)\Delta t)^\nu]
\mathbb E_{\rho} [(\sigma y)^{\alpha+\beta-\nu}]
\label{entropy:prod:add2}
\end{aligned}
\end{equation}
where $\bar{\mu}$ is the equilibrium measure for $x_i$ while $\rho$ is the Gaussian measure
of $\Delta W_i$.
Using the Isserlis-Wick formula we can compute the higher moments of multivariate Gaussian random variable
from the second-order moments. Indeed, we have
\begin{equation}
\mathbb E[y^\nu] = \mathbb E[y_1^{\nu_1}...y_d^{\nu_d}] = \mathbb E[z_1z_2...z_{|\nu|}] = \left\{\begin{matrix}
0 & \ \ \text{if}\ \ \ |\nu| \ \ \ \text{odd} \\
\sum\prod \mathbb E[z_iz_j] &\ \ \ \text{if}\ \ \ |\nu| \ \ \ \text{even}
\end{matrix} \right.
\label{Isserlis:Wick:formula}
\end{equation}
where $\sum\prod$ means summing over all distinct ways of partitioning $z_1,...,z_{|\nu|}$ into pairs.
Moreover, $\mathbb E[z_iz_j] = \Sigma_{ij} \Delta t$, hence, applying (\ref{Isserlis:Wick:formula})
into (\ref{entropy:prod:add2}) and changing the multi-index notation to the usual notation,
the entropy production rate becomes
\begin{equation}
\begin{aligned}
&EP(\Delta t) = \frac{2}{\Delta t} \sum_{k_1=1}^d\sum_{k_2=1}^d\sum_{k_3=1}^d C_{k_1k_2k_3} \left\{
\mathbb E_{\bar{\mu}} [\frac{\partial^3 V}{\partial x_{k_1}\partial x_{k_2}\partial x_{k_3}}(-\frac{1}{2}\Sigma\nabla V)_{k_1}]\Sigma_{k_2k_3} \Delta t^2 \right. \\
&\left. + \mathbb E_{\bar{\mu}} [\frac{\partial^3 V}{\partial x_{k_1}\partial x_{k_2}\partial x_{k_3}}(-\frac{1}{2}\Sigma\nabla V)_{k_2}]\Sigma_{k_1k_3} \Delta t^2
+ \mathbb E_{\bar{\mu}} [\frac{\partial^3 V}{\partial x_{k_1}\partial x_{k_2}\partial x_{k_3}}(-\frac{1}{2}\Sigma\nabla V)_{k_3}]\Sigma_{k_1k_2} \Delta t^2
+ O(\Delta t^3)\right\} \\
&+ \frac{1}{\Delta t}\sum_{k_1=1}^d\sum_{k_2=1}^d\sum_{k_3=1}^d \sum_{k_4=1}^d C_{k_1k_2k_3} \left\{
\mathbb E_{\bar{\mu}} [\frac{\partial^4 V}{\partial x_{k_1}...\partial x_{k_4}}] [\Sigma_{k_1k_2}\Sigma_{k_3k_4}
+\Sigma_{k_1k_3}\Sigma_{k_2k_4}+\Sigma_{k_1k_4}\Sigma_{k_2k_3}]\Delta t^2 + O(\Delta t^3)\right\} \\
&+ \frac{1}{\Delta t} O(\Delta t^3)\, .
\label{entropy:prod:add3}
\end{aligned}
\end{equation}
Using that $(-\frac{1}{2}\Sigma\nabla V)_{k_i} = -\frac{1}{2}\sum_{k_4=1}^d \Sigma_{k_ik_4}\frac{\partial V}{\partial x_{k_4}}$,
entropy production is rewritten as
\begin{equation}
\begin{aligned}
&EP(\Delta t) = \sum_{k_1=1}^d\sum_{k_2=1}^d\sum_{k_3=1}^d \sum_{k_4=1}^d C_{k_1k_2k_3} \left\{
\Sigma_{k_1k_2}\Sigma_{k_3k_4} \left(-\mathbb E_{\bar{\mu}} [\frac{\partial^3 V}{\partial x_{k_1}\partial x_{k_3}\partial x_{k_4}}\frac{\partial V}{\partial x_{k_2}}]
+ \mathbb E_{\bar{\mu}} [\frac{\partial^4 V}{\partial x_{k_1}...\partial x_{k_4}}]\right) \right. \\
&+ \Sigma_{k_1k_3}\Sigma_{k_2k_4} \left(-\mathbb E_{\bar{\mu}} [\frac{\partial^3 V}{\partial x_{k_1}\partial x_{k_2}\partial x_{k_4}}\frac{\partial V}{\partial x_{k_3}}]
+ \mathbb E_{\bar{\mu}} [\frac{\partial^4 V}{\partial x_{k_1}...\partial x_{k_4}}]\right) \\
&\left.+ \Sigma_{k_1k_4}\Sigma_{k_2k_3} \left(-\mathbb E_{\bar{\mu}} [\frac{\partial^3 V}{\partial x_{k_1}\partial x_{k_2}\partial x_{k_3}}\frac{\partial V}{\partial x_{k_4}}]
+ \mathbb E_{\bar{\mu}} [\frac{\partial^4 V}{\partial x_{k_1}...\partial x_{k_4}}]\right) \right\} \Delta t 
+ O(\Delta t^2)\, .
\label{entropy:prod:add4}
\end{aligned}
\end{equation}
By a simple integration by parts, we observe that for any combination $k_1,...,k_4=1,...,d$
\begin{equation}
\mathbb E_{\mu} [\frac{\partial^3 V}{\partial x_{k_1}\partial x_{k_2}\partial x_{k_3}}\frac{\partial V}{\partial x_{k_4}}]
=\mathbb E_{\mu} [\frac{\partial^4 V}{\partial x_{k_1}...\partial x_{k_4}}]
\end{equation}
where the expectation is taken with respect of $\mu$ which is the invariant measure of the
continuous-time process. However, in (\ref{entropy:prod:add4}) the expectation
is w.r.t. the invariant measure of the discrete-time process (i.e., $\bar{\mu}$
instead of $\mu$). Nevertheless, Assumption~\ref{basic:assumption}
guarantees that the alternation of the measure from $\mu$ to $\bar{\mu}$ costs an
error of order $O(\Delta t)$.
Hence, for any coefficient in (\ref{entropy:prod:add4}), we obtain that
\begin{equation}
\left|\mathbb E_{\bar{\mu}} [\frac{\partial^3 V}{\partial x_{k_1}\partial x_{k_2}\partial x_{k_3}}\frac{\partial V}{\partial x_{k_4}}]
-\mathbb E_{\bar{\mu}} [\frac{\partial^4 V}{\partial x_{k_1}...\partial x_{k_4}}]\right| \leq 2K\Delta t
\end{equation}
since the potential $V$ as well its derivatives are sufficiently smooth.
Hence, we overall showed that
\begin{equation}
EP(\Delta t) = O(\Delta t^2)
\end{equation}
which completes the proof.
\end{proof}

\begin{rmrk}
Depending on the potential function the entropy production could be
even smaller. For instance, when the potential $V$ is a quadratic function (i.e. the
continuous-time process is an Ornstein-Uhlenbeck process), then, it is easily
checked by a trivial calculation of (\ref{action:func:add}) that the GC action
function is a boundary term, thus, the entropy production of the explicit EM
scheme is zero. However, for a generic potential $V$ we expect that  
the entropy production rate decays quadratically as a function of
$\Delta t$ but not faster.
\end{rmrk}

\subsubsection{Fourth-order potential on a torus}
\begin{figure}[!htb]
\begin{center}
\includegraphics[width=.8\textwidth]{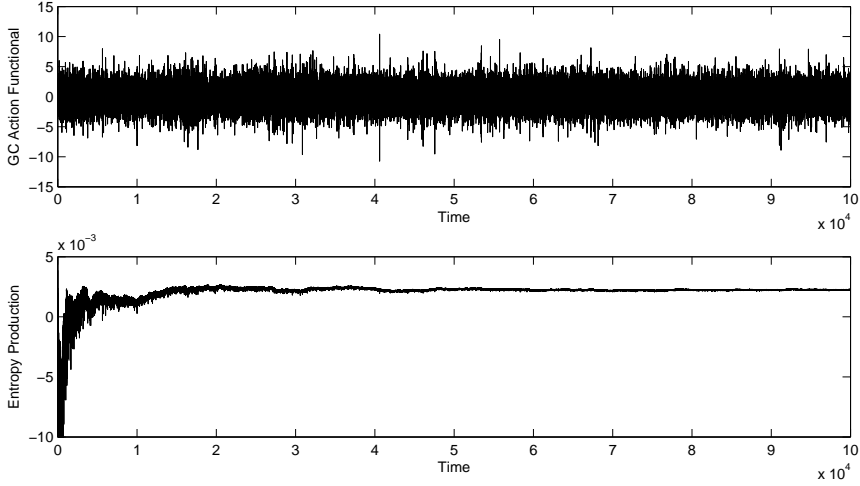}
\caption{Upper panel: The GC action functional as a function of time for fixed $\Delta t=0.05$.
Its variance is large necessitating the use of many samples in order to obtain statistically
confident quantities. Lower Panel: The entropy production rate as a function of time for the same
$\Delta t$. It converges to a positive value as expected.}
\label{gcaf:ep:fig}
\end{center}
\end{figure}

\begin{figure}[!htb]
\begin{center}
\includegraphics[width=.8\textwidth]{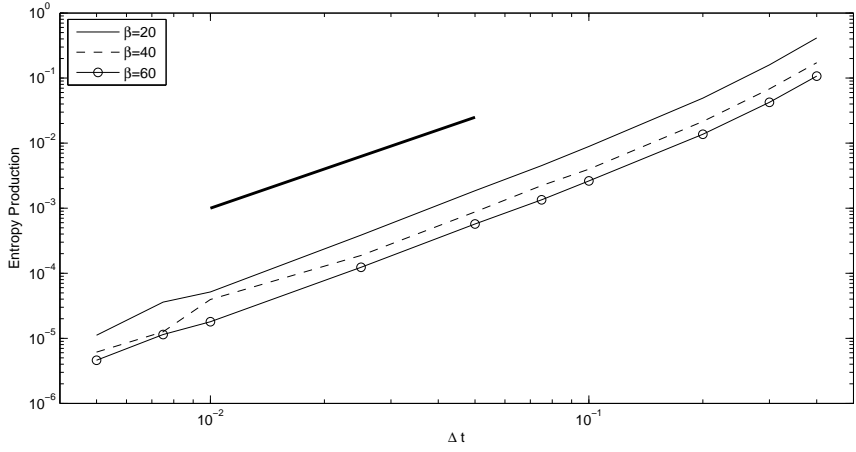}
\caption{Entropy production rate as a function of time step $\Delta t$ for additive noise.
The entropy production rate is of order $O(\Delta t^2)$ for small $\Delta t$ while it
decreases linearly as a function of inverse temperature $\beta$.}
\label{add:fig}
\end{center}
\end{figure}

Lets now proceed with an important example where the potential is a forth-order
polynomial while the process takes values on a torus. Assume $d=2$ while potential
$V=V_\beta$ is given by
\begin{equation}
V_\beta(x) = \beta\left(\frac{|x|^4}{4} - \frac{|x|^2}{2}\right)
\end{equation}
where $\beta$ is a positive real number which in statistical mechanics has the meaning of
the inverse temperature. The diffusion matrix is set to $\sigma=\sqrt{2\beta^{-1}}I_d$.
Based on \cite{Mattingly:02}, Assumption~\ref{basic:assumption} is satisfied
because the domain is restricted to a torus, the potential is locally Lipschitz
continuous and the covariance matrix is elliptic.
Figure~\ref{gcaf:ep:fig} presents both the GC action functional (upper panel)
and the entropy production rate (lower panel) as a function of time for fixed $\Delta t=0.05$.
Both quantities are numerically computed while the inverse temperature is set
to $\beta=10$. Even though the variance of the GC action functional is large, entropy
production which is the cumulative sum of the GC functional converges due to the
law of large numbers to a (positive) value after relatively long time. Additionally,
due to the ergodicity assumption, it converges to the correct value.
Figure~\ref{add:fig} shows the loglog plot of the numerical entropy production
rate as a function of $\Delta t$ for $\beta=20,\ 40,\ 60$. Final time was set to
$t=2\cdot10^6$ while initial point was set to one of the attraction points of the
deterministic counterpart. For reader's convenience, the thick black line denotes the
$O(\Delta t^2)$ rate of convergence. This plot is in agreement with the theorem's
estimate (\ref{entr:prod:add}).
Notice also that, for small $\Delta t$, entropy production rate is very close to 0 and
even larger final time is needed in order to obtain a statistically confident numerical
estimate for the entropy production. Moreover, as it is evident from the figure and
the GC action functional in (\ref{action:func:add}), the dependence of the entropy
production w.r.t. the inverse temperature is inverse proportional. Thus, from a
statistical mechanics point of view, the larger is the temperature the larger
--in a linear manner-- is the entropy production rate of the numerical scheme.

\subsection{Entropy Production for the Multiplicative Noise Case: Euler-Marayuma scheme}\label{multiplicative}
In this section we consider the EM scheme for  overdamped Langevin processes with multiplicative noise.
For simplicity  we restrict our discussion to the one dimensional case,  but our results extend immediately to higher dimension if the the diffusion matrix $\sigma(x)$ is diagonal.
%
%
%Nonetheless, the results and conclusions of this subsection
%for both explicit EM and Milstein's schemes are valid in a more general, multi-dimensional  setting
%where the diffusion matrix $\sigma(x)$ is diagonal.
%
%In order to study the entropy production rate of the explicit EM scheme for the overdamped Langevin
%process with multiplicative noise,
% the remainder terms of the GC action functional should be studied.
%In this direction 
We rewrite the GC action function given in Lemma~\ref{overdamped:Langevin:lemma}, % for $1d$
\begin{equation}
\begin{aligned}
W(n;\Delta t) &\dot{=} - \frac{1}{2}\sum_{i=0}^{n-1} [V'(x_{i+1})+V'(x_i)]\Delta x_{i}
+ \frac{1}{2}  \sum_{i=0}^{n-1} [\Sigma^{-1}(x_{i+1})\Sigma'(x_{i+1}) + \Sigma^{-1}(x_i)\Sigma'(x_i)]\Delta x_{i} \\
&\quad +\frac{1}{2\Delta t} \sum_{i=0}^{n-1} \left[\Sigma^{-1}(x_{i+1})-\Sigma^{-1}(x_i)\right]\Delta x_i^2 \\
&=: W_1(n;\Delta t) + W_2(n;\Delta t) + W_3(n;\Delta t)\, .
\end{aligned}
\end{equation}
The first term $W_1(n;\Delta t)$ has been computed in the previous section and after an
interesting and rather unexpected cancellation it was proved to be of order $O(\Delta t^2)$.
For the multiplicative case, a cancellation also occurs 
(see (\ref{cancel:multi}) and (\ref{remain:term:multi}) below) but it does not fully eliminate the lower order term; 
in the end $W_1(n;\Delta t)$ contributes to the entropy
production an  $O(\Delta t)$ term. Additionally, $W_2(n;\Delta t)$ turns out to be  the sum of gradient terms
since  %$\Sigma(x)\in\mathbb R$ and  
$\Sigma^{-1}(x)\Sigma'(x) = (\log\Sigma(x))'$.
Thus, assuming a suitable condition on $\Sigma(x)$, the same computation as for $W_1(n;\Delta t)$
applies and the entropy production rate stemming from $W_2(n;\Delta t)$ is also of order $O(\Delta t)$.
However, $W_3(n;\Delta t)$ contributes to the entropy production  a nonzero, positive term which is of order $O(1)$.
The following theorem summarizes the behavior of entropy production rate for the explicit
EM scheme for multiplicative noise.
\begin{thrm}
Let Assumption~\ref{basic:assumption} hold and assume  that the potential
function $V$ has a bounded fifth-order derivative,  while there exists $M>0$ such that
$\Sigma(x)>M^{-1}$ for all $x$. \\
(a) If  $c:=\frac{3}{4}\mathbb E_{\mu} [(\Sigma^{-1})(x)(\Sigma')^2(x)]$, then, for sufficiently
small $\Delta t$, there exists $C=C(V,\Sigma)>0$ independent of $\Delta t$ such that
\begin{equation}
|EP(\Delta t) - c| \leq C \Delta t
\label{entr:prod:multi_a}
\end{equation}
(b) Assuming that $\mathbb E_{\mu} [(\Sigma^{-1})(x)(\Sigma')^2(x)]\neq0$, then, for sufficiently
small $\Delta t$, there exists a lower bound $c'=c'(V,\Sigma)>0$ independent of $\Delta t$ such that
\begin{equation}
c' \leq EP(\Delta t)
\label{entr:prod:multi}
\end{equation}
\label{entr:prod:multi:theorem}
\end{thrm}
\begin{proof}
The assumption  that $\Sigma(x)>M^{-1}\ \forall x$, which is the ellipticity condition  in one space dimension,
is necessary because it implies that  $\Sigma^{-1}(x)$ as well its derivatives are bounded around 0.
Additionally, as discussed earlier both $W_1(n;\Delta t)$ and $W_2(n;\Delta t)$ contribute to the entropy
production by a $O(\Delta t)$ amount. Therefore we can concentrate on the term  $W_3(n;\Delta t)$;  after a Taylor series
expansion we have,
\begin{equation*}
\begin{aligned}
&W_3(n;\Delta t) = \frac{1}{2\Delta t} \sum_{i=0}^{n-1} \left[(\Sigma^{-1})'(x_i)\Delta x_i^3 + \frac{1}{2} (\Sigma^{-1})''(x_i)\Delta x_i^4
+ \frac{1}{2\Delta t} \sum_{i=0}^{n-1} \int_0^1(1-t)(\Sigma^{-1})'''(tx_{i+1}+(1-t)x_i)dt \Delta x_i^5 \right]\\
&= \frac{1}{2\Delta t} \sum_{i=0}^{n-1}\sum_{k=0}^3 \binom{3}{k} (\Sigma^{-1})'(x_i)
(-\frac{1}{2}\Sigma(x_i) V'(x_i)\Delta t + \frac{1}{2} \Sigma'(x_i)\Delta t)^k (\sigma(x_i) \Delta W_i)^{3-k} \\
&+ \frac{1}{4\Delta t} \sum_{i=0}^{n-1}\sum_{k=0}^4 \binom{4}{k} (\Sigma^{-1})''(x_i)
(-\frac{1}{2}\Sigma(x_i) V'(x_i)\Delta t + \frac{1}{2} \Sigma'(x_i)\Delta t)^k (\sigma(x_i) \Delta W_i)^{4-k} \\
&+ \frac{1}{2\Delta t} \sum_{i=0}^{n-1}\sum_{k=0}^5 \binom{5}{k}\int_0^1(1-t)(\Sigma^{-1})'''(tx_{i+1}+(1-t)x_i)dt
(-\frac{1}{2}\Sigma(x_i) V'(x_i)\Delta t + \frac{1}{2} \Sigma'(x_i)\Delta t)^k (\sigma(x_i) \Delta W_i)^{5-k} \, .
\end{aligned}
\end{equation*}
As in Theorem~\ref{add:noise:theorem}, applying the ergodic lemmas of the appendix,
the entropy production rate stemming from $W_3(n;\Delta t)$ equals to
\begin{equation}
\begin{aligned}
&EP_3(\Delta t) = \lim_{t\rightarrow\infty} \frac{W_3(n;\Delta t)}{n\Delta t} \\
&= \frac{1}{2\Delta t^2} \sum_{k=0}^3 \binom{3}{k} \mathbb E_{\bar{\mu}} [(\Sigma^{-1})'(x)
(-\frac{1}{2}\Sigma(x) V'(x)\Delta t + \frac{1}{2} \Sigma'(x)\Delta t)^k \sigma(x)^{3-k}] \mathbb E_\rho [\Delta W^{3-k}] \\
&+ \frac{1}{4\Delta t^2} \sum_{k=0}^4 \binom{4}{k} \mathbb E_{\bar{\mu}} [(\Sigma^{-1})''(x)
(-\frac{1}{2}\Sigma(x) V'(x)\Delta t + \frac{1}{2} \Sigma'(x)\Delta t)^k \sigma(x)^{4-k}] \mathbb E_\rho [\Delta W^{4-k}] \\
&+ \frac{1}{2\Delta t^2} \sum_{k=0}^5 \mathbb E_{\bar{\mu}\times\rho} [R(x,y)
(-\frac{1}{2}\Sigma(x) V'(x)\Delta t + \frac{1}{2} \Sigma'(x)\Delta t)^k \sigma(x)^{5-k}] \mathbb E_\rho [\Delta W^{5-k}] \\
&= \frac{1}{2\Delta t^2} \left[-\frac{3}{2}\mathbb E_{\bar{\mu}} [(\Sigma^{-1})'(x)\Sigma^2(x) V'(x)]\Delta t^2
+ \frac{3}{2} \mathbb E_{\bar{\mu}} [(\Sigma^{-1})'(x)\Sigma'(x)\Sigma(x)]\Delta t^2 + O(\Delta t^3)\right] \\
&+ \frac{1}{4\Delta t^2} \left[\mathbb E_{\bar{\mu}} [(\Sigma^{-1})''(x)\Sigma^2(x)]3\Delta t^2 + O(\Delta t^3)\right]
+ \frac{1}{2\Delta t^2} O(\Delta t^3) \\
&= \frac{3}{4}\left[-\mathbb E_{\bar{\mu}} [(\Sigma^{-1})'(x)\Sigma^2(x) V'(x)] + \frac{1}{2}\mathbb E_{\bar{\mu}} [(\Sigma^{-1})'(x)(\Sigma^2)'(x)]
+ \mathbb E_{\bar{\mu}} [(\Sigma^{-1})''(x)\Sigma^2(x)] \right] + O(\Delta t)
\end{aligned}
\end{equation}
On the other hand, % it holds for the invariant measure $\mu$ that
\begin{equation}
\mathbb E_{\mu} [(\Sigma^{-1})'(x)\Sigma^2(x) V'(x)] = \mathbb E_{\mu} [(\Sigma^{-1})''(x)\Sigma^2(x)]
+ \mathbb E_{\mu} [(\Sigma^{-1})'(x)(\Sigma^2)'(x)]
\label{cancel:multi}
\end{equation}
Using  (\ref{stat:average:estimate}) in Assumption~\ref{basic:assumption} we obtain, as in
the additive case, that
\begin{equation}
%\begin{aligned}
%&EP_3(\Delta t) = -\frac{3}{8} \mathbb E_{\bar{\mu}} [(\Sigma^{-1})'(x)(\Sigma^2)'(x)] + O(\Delta t) \\
%\Rightarrow &
EP_3(\Delta t)- \frac{3}{4} \mathbb E_{\bar{\mu}} [(\Sigma^{-1})(x)(\Sigma')^2(x)] = O(\Delta t)
\label{remain:term:multi}
%\end{aligned}
\end{equation}
which concludes the proof of (a). Part (b) is a direct consequence of (a).
\end{proof}

\subsubsection{Example: Quadratic potential on $\mathbb R$}
\begin{figure}[!htb]
\begin{center}
\includegraphics[width=.8\textwidth]{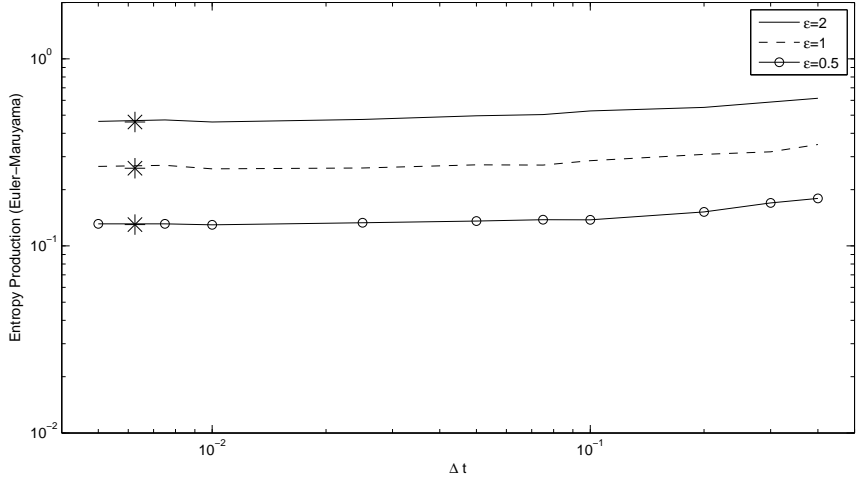}
\caption{Entropy production rate as a function of time step $\Delta t$ for multiplicative noise and
the explicit EM scheme. As Theorem~\ref{entr:prod:multi:theorem} asserts, entropy production
does not decrease as $\Delta t$ is decreased. This results in a permanent loss of reversibility
which cannot be fixed by reducing the time step. Star symbols denote the theoretical value of
the lower bound as it is given by the Theorem
(i.e., $c'\approx c=\frac{3}{4}\mathbb E_\mu[(\Sigma_\epsilon)^{-1}(x) (\Sigma_\epsilon')^2(x)]$). The
agreement between the theoretical and the numerical values is excellent.}
\label{multi:fig1}
\end{center}
\end{figure}

Let the quadratic potential  $V(x)=\frac{x^2}{2}$ , and the   diffusion term 
\begin{equation}
\sigma_\epsilon(x) = \sqrt{\frac{1}{1+\epsilon x^2}}\, .
\end{equation}
The choice of the diffusion term is justified by the fact that we can control its variation in terms
of $x$, while  sending $\epsilon$ to zero, the additive noise case is recovered.
The invariant measure of this process is the Gaussian measure with zero mean and variance one.
%This invariant measure is the simplest measure to be considered.
Moreover, all the assumptions of Theorem~\ref{entr:prod:multi:theorem} are satisfied
thus we expect a $O(1)$ behavior of the entropy production rate at least for small $\Delta t$.
Indeed, Figure~\ref{multi:fig1} shows the  numerically-computed entropy
production as a function of $\Delta t$, which clearly  does not decrease to
zero as $\Delta t$ tends to zero. Consequently, the explicit EM scheme for the multiplicative noise case totally
destroys the reversibility property of the discrete-time approximation process independently
of how small time-step is selected.
Additionally, notice that as $\epsilon$ decreases, entropy production also decreases. This behavior is
 expected since $\sigma(x)\rightarrow\sigma=\text{constant}$ as $\epsilon\rightarrow0$
and in combination with the quadratic potential $V$, $EP(\Delta t)\rightarrow0$ as
$\epsilon\rightarrow0$ for any $\Delta t$ sufficiently small.

\subsection{ Entropy Production for the Multiplicative Noise Case: Milstein scheme}
Since the EM scheme has entropy production rate which does not decrease as $\Delta t$ decreases,
we turn our attention to  %an immediate question to ask is what happens when 
%a higher-order scheme such as  
the  Milstein's scheme which  is the next higher-order scheme \cite{Kloeden:99, Milstein:04}:% and its explicit version is
%given by
\begin{equation}
\begin{aligned}
&x_{i+1} = x_i -\frac{1}{2}\Sigma(x_i) V'(x_i)\Delta t + \frac{1}{2} \Sigma'(x_i)\Delta t + \sigma(x_i) \Delta W_i
+ \frac{1}{2}\sigma(x_i)\sigma'(x_i)(\Delta W_i^2-\Delta t)\, , \\
\label{multi:sde:explMilstein}
\end{aligned}
\end{equation}
which can be  rewritten as
\begin{equation}
\Delta x_i = a(x_i)\Delta t + \sigma(x_i) \Delta W_i + \frac{1}{4}\Sigma'(x_i) \Delta W_i^2\, ,
\end{equation}
where $a(x_i)=-\frac{1}{2}\Sigma(x_i) V'(x_i) + \frac{1}{4} \Sigma'(x_i)$.
Since $\Delta W_i$ is a zero-mean Gaussian random variable with variance $\Delta t$, the transition
probability for Milstein's scheme is
\begin{equation}
\begin{aligned}
\Pi(x_i,x_{i+1}) &= \frac{1}{|\sqrt{2\pi\Delta t Z(x_i,\Delta x_i)}|}
\left[\exp\left(-\frac{1}{2\Delta t}\left|\frac{-\sigma(x_i) + \sqrt{Z(x_i,\Delta x_i)}}{\frac{1}{2}\Sigma'(x_i)}\right|^2\right)\right. \\
&\left.+ \exp\left(-\frac{1}{2\Delta t}\left|\frac{\sigma(x_i) + \sqrt{Z(x_i,\Delta x_i)}}{\frac{1}{2}\Sigma'(x_i)}\right|^2\right)\right]
\label{trans:prob:Milstein}
\end{aligned}
\end{equation}
where
\begin{equation}
Z(x_i,\Delta x_i) =  \Sigma(x_i) + \Sigma'(x_i)\left(\Delta x_i - a(x_i)\Delta t\right)\, .
\label{norm:factor:Milstein}
\end{equation}
Notice also that $Z(x_i,\Delta x_i) = (\sigma(x_i)+\frac{1}{2}\Sigma'(x_i)\Delta W_i)^2\geq0$
which is positive almost surely. Moreover, the arguments of the exponentials in (\ref{trans:prob:Milstein})
are of different order in terms of $\Delta t$. Indeed, it is straightforward to show that for small time
step, $\Delta t$, the argument of the first exponential in (\ref{trans:prob:Milstein}) is of order $O(1)$
while the argument of the second exponential is of order $O(\frac{1}{\Delta t})$. Thus, as $\Delta t$
tends to zero, the second exponential becomes exponentially small and the dominating term is
the first exponential. Therefore, using the fact that $\log\left(e^{-a}+e^{-b/\Delta t}\right)=-a+O(e^{-b/\Delta t})$
for positive $a$ and $b$, the GC action functional for Milstein's scheme reduces  to
\begin{equation}
\begin{aligned}
W(n;\Delta t) &= -\frac{1}{2}\sum_{i=0}^{n-1} \log \frac{Z(x_{i},\Delta x_i)}{Z(x_{i+1},-\Delta x_i)}
-\frac{2}{\Delta t} \sum_{i=0}^{n-1} \left[\left(\frac{-\sigma(x_i) + \sqrt{Z(x_i,\Delta x_i)}}{\frac{1}{2}\Sigma'(x_i)}\right)^2
-\left(\frac{-\sigma(x_{i+1}) + \sqrt{Z(x_{i+1},-\Delta x_i)}}{\frac{1}{2}\Sigma'(x_{i+1})}\right)^2 \right] \\
&= W_1(n;\Delta t) + W_2(n;\Delta t)
\label{GC:action:func:Milstein}
\end{aligned}
\end{equation}
where $Z(x_{i+1},-\Delta x_i) = \Sigma(x_{i+1}) + \Sigma'(x_{i+1})\left(-\Delta x_i - a(x_{i+1})\Delta t\right)$.
The following theorem demonstrates  that the entropy production of the Milstein Scheme is at least
linear in $\Delta t$:

\begin{thrm}
Under the assumptions of Theorem~\ref{entr:prod:multi:theorem} and for sufficiently
small $\Delta t$, there exists $C=C(V,\Sigma)>0$ independent of $\Delta t$ such that
\begin{equation}
EP(\Delta t) \leq C\Delta t
\label{entr:prod:multi:Milstein}
\end{equation}
\label{Milstein:entr:prod:multi:theorem}
\end{thrm}
\begin{proof}
In order to compute the detailed asymptotics for  $W_1(n;\Delta t)$ and
$W_2(n;\Delta t)$ we write the partition function $Z(x_i, \Delta x_i)$ as 
\begin{equation}
\begin{aligned}
Z(x_i, \Delta x_i) &= \Sigma(x_i) + \Sigma'(x_i)\left(\Delta x_i - a(x_i)\Delta t\right) \\
&= \Sigma(x_{i+1}) - \left( \frac{1}{2}\Sigma''(x_i)\Delta x_i^2 + \frac{1}{6}\Sigma'''(x_i)\Delta x_i^3 + \Sigma'(x_i)a(x_i)\Delta t \right) + O(\Delta x_i^4) \ .
\end{aligned}
\end{equation}
Similarly we have
\begin{equation}
\begin{aligned}
Z(x_{i+1}, -\Delta x_i) &= \Sigma(x_{i+1}) + \Sigma'(x_{i+1})\left(-\Delta x_i - a(x_{i+1})\Delta t \right) \\
&= \Sigma(x_i) - \left( \frac{1}{2}\Sigma''(x_{i+1})\Delta x_i^2 - \frac{1}{6}\Sigma'''(x_{i+1})\Delta x_i^3 + \Sigma'(x_{i+1})a(x_{i+1})\Delta t \right) + O(\Delta x_i^4) \ ,
\end{aligned}
\end{equation}
and thus 
\begin{equation}
\begin{aligned}
&Z(x_{i+1}, -\Delta x_i) - Z(x_{i-1}, \Delta x_{i-1}) = -\frac{1}{2} (\Sigma''(x_{i+1})\Delta x_i^2 - \Sigma''(x_{i-1})\Delta x_{i-1}^2) \\
&+ \frac{1}{6}(\Sigma'''(x_{i+1})\Delta x_i^3 + \Sigma'''(x_{i-1})\Delta x_{i-1}^3) - (\Sigma'(x_{i+1})a(x_{i+1})-\Sigma'(x_{i-1})a(x_{i-1}))\Delta t
\end{aligned}
\end{equation}
is obtained. Moreover, in what follows and by slight abuse of $O(\cdot)$ notation,
we repeatedly use the relation
\begin{equation}
\left[f(x_i)g(x_{i\pm1}) - f(x_{i\pm1})g(x_i)\right]\Delta x_i^k = O(\Delta x_i^{k+1})
\label{simplification:eq}
\end{equation}
which holds for any $i,k=0,1,...$ and any smooth functions $f$ and $g$ and it is
easily derived by suitable Taylor expansions of the functions. We obtain for 
$W_1(n;\Delta t)$ 
\begin{equation}
\begin{aligned}
W_1(n;\Delta t) &= \frac{1}{2}\sum_{i=0}^{n-1} \log\frac{Z(x_{i+1}, -\Delta x_i)}{Z(x_i, \Delta x_i)} \\
&\dot{=} \frac{1}{2}\sum_{i=0}^{n-1} \log\frac{Z(x_{i+1}, -\Delta x_i)}{Z(x_{i-1}, \Delta x_{i-1})} \\
&= \frac{1}{2}\sum_{i=0}^{n-1} \log\left(1-\frac{\frac{1}{2} (\Sigma''(x_{i+1})\Delta x_i^2 - \Sigma''(x_{i-1})\Delta x_{i-1}^2)
+(\Sigma'(x_{i+1})a(x_{i+1})-\Sigma'(x_{i-1})a(x_{i-1}))\Delta t + O(\Delta x_i^3)}{Z(x_{i-1}, \Delta x_{i-1})}   \right) \\
&= \frac{1}{2}\sum_{i=0}^{n-1} \sum_{k=1}^\infty \left( \frac{\frac{1}{2} (\Sigma''(x_{i+1})\Delta x_i^2 - \Sigma''(x_{i-1})\Delta x_{i-1}^2)
+ O(\Delta t\Delta x_i + \Delta x_i^3)}{Z(x_{i-1}, \Delta x_{i-1})}   \right)^k \\
&= \frac{1}{4}\sum_{i=0}^{n-1} \left[\frac{\Sigma''(x_{i+1})\Delta x_i^2 - \Sigma''(x_{i-1})\Delta x_{i-1}^2}
{\Sigma(x_i) + O(\Delta t + \Delta x_i^2)} + O(\Delta t\Delta x_i + \Delta x_i^3) \right] \\
&= \frac{1}{4}\sum_{i=0}^{n-1} \left[\frac{\Sigma''(x_{i+1})\Delta x_i^2 - \Sigma''(x_{i-1})\Delta x_{i-1}^2}{\Sigma(x_i)}
\sum_{k=0}^\infty \left( O(\Delta t + \Delta x_i^2) \right)^k + O(\Delta t\Delta x_i + \Delta x_i^3) \right] \\
&= \frac{1}{4}\sum_{i=0}^{n-1} \left[\frac{\Sigma''(x_{i+1})\Delta x_i^2 - \Sigma''(x_{i-1})\Delta x_{i-1}^2}{\Sigma(x_i)}
+ O(\Delta t\Delta x_i + \Delta x_i^3) \right] \\
&\dot{=} \frac{1}{4}\sum_{i=0}^{n-1} \left[\left(\frac{1}{\Sigma(x_i)} - \frac{1}{\Sigma(x_{i+1})} \right) \Sigma''(x_{i})\Delta x_i^2
+ O(\Delta t\Delta x_i + \Delta x_i^3) \right] \\
&= \sum_{i=0}^{n-1} O(\Delta x_i^3) + \Delta t \sum_{i=0}^{n-1} O(\Delta x_i)
\end{aligned}
\end{equation}

The second term of the GC action functional is rewritten as
\begin{equation}
\begin{aligned}
W_2(n;\Delta t) &= \frac{2}{\Delta t} \sum_{i=0}^{n-1} \left[\left(\frac{\sigma(x_{i+1}) - \sqrt{Z(x_{i+1},-\Delta x_i)}}{\frac{1}{2}\Sigma'(x_{i+1})}\right)^2
- \left(\frac{\sigma(x_i) - \sqrt{Z(x_i,\Delta x_i)}}{\frac{1}{2}\Sigma'(x_i)}\right)^2 \right] \\
&= \frac{2}{\Delta t} \sum_{i=0}^{n-1} \left[\left(\frac{\sigma(x_{i+1}) - \sigma(x_i)\Big(1-\frac{1}{2\Sigma(x_i)}\big(\frac{1}{2}\Sigma''(x_{i+1})\Delta x_i^2
-\frac{1}{6}\Sigma'''(x_{i+1})\Delta x_i^3 +\Sigma'(x_{i+1})a(x_{i+1})\Delta t + O(\Delta x_i^4) \big)\Big)}{\frac{1}{2}\Sigma'(x_{i+1})}\right)^2 \right. \\
&\left. - \left(\frac{\sigma(x_i) - \sigma(x_{i+1})\Big(1- \frac{1}{2\Sigma(x_{i+1})}\big(\frac{1}{2}\Sigma''(x_i)\Delta x_i^2
+\frac{1}{6}\Sigma'''(x_i)\Delta x_i^3 +\Sigma'(x_i)a(x_i)\Delta t + O(\Delta x_i^4) \big)\Big)}{\frac{1}{2}\Sigma'(x_i)}\right)^2 \right] \\
&= \frac{2}{\Delta t} \sum_{i=0}^{n-1} \left[\left(2\frac{\sigma(x_{i+1}) - \sigma(x_i)}{\Sigma'(x_{i+1})} + \frac{1}{2}\frac{\Sigma''(x_{i+1})\Delta x_i^2}{\sigma(x_i)\Sigma'(x_{i+1})}
- \frac{1}{6} \frac{\Sigma'''(x_{i+1})\Delta x_i^3}{\sigma(x_i)\Sigma'(x_{i+1})} + \frac{1}{2}\frac{\Sigma'(x_{i+1})a(x_{i+1})\Delta t}{\sigma(x_i)\Sigma'(x_{i+1})}
+ O(\Delta x_i^4) \right)^2 \right. \\
&\left. - \left(2\frac{\sigma(x_i) - \sigma(x_{i+1})}{\Sigma'(x_i)} + \frac{1}{2}\frac{\Sigma''(x_i)\Delta x_i^2}{\sigma(x_{i+1})\Sigma'(x_i)}
+ \frac{1}{6} \frac{\Sigma'''(x_i)\Delta x_i^3}{\sigma(x_{i+1})\Sigma'(x_i)} + \frac{1}{2}\frac{\Sigma'(x_i)a(x_i)\Delta t}{\sigma(x_{i+1})\Sigma'(x_i)}
+ O(\Delta x_i^4) \right)^2 \right] \\
\end{aligned}
\end{equation}
where a Taylor series expansion to the square root function was applied.
Expanding the squares and keeping only the terms that have order in terms of $\Delta x_i$
less than 5 we obtain that
\begin{equation}
\begin{aligned}
W_2(n;\Delta t) &= \frac{2}{\Delta t} \sum_{i=0}^{n-1} \left[
4\left(\frac{(\sigma(x_{i+1}) - \sigma(x_i))^2}{\Sigma'(x_{i+1})^2} - \frac{(\sigma(x_{i+1}) - \sigma(x_i))^2}{\Sigma'(x_i)^2} \right) \right. \\
&+2\left(\frac{\Sigma''(x_{i+1})(\sigma(x_{i+1}) - \sigma(x_i))}{\sigma(x_i)\Sigma'(x_{i+1})^2}
+ \frac{\Sigma''(x_i)(\sigma(x_{i+1}) - \sigma(x_i))}{\sigma(x_{i+1})\Sigma'(x_i)^2} \right) \Delta x_i^2 \\
&-\frac{2}{3}\left(\frac{\Sigma'''(x_{i+1})(\sigma(x_{i+1}) - \sigma(x_i))}{\sigma(x_i)\Sigma'(x_{i+1})^2}
-  \frac{\Sigma'''(x_i)(\sigma(x_{i+1}) - \sigma(x_i))}{\sigma(x_{i+1})\Sigma'(x_i)^2} \right)\Delta x_i^3 \\
&+ 4\left(\frac{a(x_{i+1})(\sigma(x_{i+1}) - \sigma(x_i))}{\sigma(x_i)\Sigma'(x_{i+1})}
+ \frac{a(x_i)(\sigma(x_{i+1}) - \sigma(x_i))}{\sigma(x_{i+1})\Sigma'(x_i)} \right) \Delta t \\
&\left.+ \left(\frac{\Sigma''(x_{i+1})a(x_{i+1})}{\Sigma(x_i)\Sigma'(x_{i+1})}
- \frac{\Sigma''(x_i)a(x_i)}{\Sigma(x_{i+1})\Sigma'(x_i)} \right) \Delta t\Delta x_i^2 + O(\Delta x_i^5) + O(\Delta t\Delta x_i^3) \right] \\
&= \frac{2}{\Delta t} \sum_{i=0}^{n-1} \left[
\frac{\sigma'(x_{i+1})^2\Delta x_i^2 - \sigma'(x_{i+1})\sigma''(x_{i+1})\Delta x_i^3
+ \left(\frac{1}{3}\sigma'(x_{i+1})\sigma'''(x_{i+1})+\frac{1}{4}\sigma''(x_{i+1})^2\right)\Delta x_i^4}{\sigma(x_{i+1})^2\sigma'(x_{i+1})^2} \right. \\
&- \frac{\sigma'(x_i)^2\Delta x_i^2 + \sigma'(x_i)\sigma''(x_i)\Delta x_i^3
+ \left(\frac{1}{3}\sigma'(x_i)\sigma'''(x_i)+\frac{1}{4}\sigma''(x_i)^2\right)\Delta x_i^4}{\sigma(x_i)^2\sigma'(x_i)^2} \\
&+2\frac{\sigma(x_{i+1})(\sigma'(x_{i+1})\Delta x_i - \frac{1}{2}\sigma''(x_{i+1})\Delta x_i^2)\Sigma''(x_{i+1})\Sigma'(x_i)^2
+ \sigma(x_i)(\sigma'(x_i)\Delta x_i + \frac{1}{2}\sigma''(x_i)\Delta x_i^2)\Sigma''(x_i)\Sigma'(x_{i+1})^2}
{\sigma(x_i)\sigma(x_{i+1})\Sigma'(x_i)^2\Sigma'(x_{i+1})^2}\Delta x_i^2 \\
&\left. + 2\left(\frac{a(x_{i+1})(\sigma'(x_{i+1})\Delta x_i - \frac{1}{2}\sigma''(x_{i+1})\Delta x_i^2)}{\sigma(x_{i+1})\sigma'(x_{i+1})\sigma(x_i)}
+ \frac{a(x_i)(\sigma'(x_i)\Delta x_i + \frac{1}{2}\sigma''(x_i)\Delta x_i^2)}{\sigma(x_i)\sigma'(x_i)\sigma(x_{i+1})}
\right)\Delta t + O(\Delta x_i^5) + O(\Delta t\Delta x_i^3) \right] \\
&= \frac{2}{\Delta t} \sum_{i=0}^{n-1} \left[
\left(\frac{1}{\sigma(x_{i+1})^2} - \frac{1}{\sigma(x_i)^2}\right)\Delta x_i^2
- \left(\frac{\sigma''(x_{i+1})}{\sigma(x_{i+1})^2\sigma'(x_{i+1})} + \frac{\sigma''(x_i)}{\sigma(x_i)^2\sigma'(x_i)} \right)\Delta x_i^3 \right. \\
&\left. + \frac{\Sigma''(x_{i+1})\Sigma'(x_i) + \Sigma''(x_i)\Sigma'(x_{i+1})}{\sigma(x_i)\sigma(x_{i+1})\Sigma'(x_i)\Sigma'(x_{i+1})}\Delta x_i^3
+ 2\frac{a(x_{i+1})+a(x_i)}{\sigma(x_i)\sigma(x_{i+1})}\Delta x_i\Delta t + O(\Delta x_i^5) + O(\Delta t\Delta x_i^3) \right] \\
&= \frac{2}{\Delta t} \sum_{i=0}^{n-1} \left[-\left(\frac{\sigma'(x_{i+1})}{\sigma(x_{i+1})^3} + \frac{\sigma'(x_i)}{\sigma(x_i)^3}
+ \frac{\sigma''(x_{i+1})}{\sigma(x_{i+1})^2\sigma'(x_{i+1})} + \frac{\sigma''(x_i)}{\sigma(x_i)^2\sigma'(x_i)} \right)\Delta x_i^3 \right. \\
&\left. + \left(\frac{\sigma'(x_{i+1})}{\sigma(x_{i+1})^2\sigma(x_i)} + \frac{\sigma''(x_{i+1})}{\sigma(x_{i+1})\sigma'(x_{i+1})\sigma(x_i)}
+ \frac{\sigma'(x_i)}{\sigma(x_i)^2\sigma(x_{i+1})} + \frac{\sigma''(x_i)}{\sigma(x_i)\sigma'(x_i)\sigma(x_{i+1})} \right)\Delta x_i^3 + O(\Delta x_i^5) \right] \\
&+ 4 \sum_{i=0}^{n-1} \left[ \left(\frac{a(x_{i+1})}{\Sigma(x_{i+1})} + \frac{a(x_{i})}{\Sigma(x_{i})} \right)\Delta x_i + O(\Delta x_i^3) \right]
\end{aligned}
\end{equation}

After few more Taylor expansions in the first sum, the terms of order $\Delta x_i^3$ are cancelled
out while the forth order terms consists of differences of the form (\ref{simplification:eq}) thus they
become fifth order. Moreover, the second sum can be handled exactly as the terms $W_1$ and $W_2$ in 
EM scheme using (\ref{gen:trapez:rule}) and the leading term
is of order $O(\Delta x_i^3)$. Overall, we rigorously computed that
\begin{equation}
W_2(n;\Delta t) \dot{=} \frac{1}{\Delta t} \sum_{i=0}^{n-1} O(\Delta x_i^5) + \sum_{i=0}^{n-1} O(\Delta x_i^3) \ .
\end{equation}

Therefore, the entropy production for Milstein's scheme in the one dimensional overdamped Langevin
case with multiplicative noise is at least of order
\begin{equation}
\begin{aligned}
EP(\Delta t) &= \lim_{t\rightarrow\infty} \frac{1}{n\Delta t} (W_1(n;\Delta t) + W_2(n;\Delta t)) \\
&= \lim_{n\rightarrow\infty} \frac{1}{n}\sum_{i=0}^{n-1} O(\Delta x_i)
+ \frac{1}{\Delta t} \lim_{n\rightarrow\infty} \frac{1}{n}\sum_{i=0}^{n-1} O(\Delta x_i^3)
+ \frac{1}{\Delta t^2} \lim_{n\rightarrow\infty} \frac{1}{n}\sum_{i=0}^{n-1} O(\Delta x_i^5) \\
&= O(\Delta t) + \frac{1}{\Delta t} O(\Delta t^2) + \frac{1}{\Delta t^2} O(\Delta t^3) \\
&= O(\Delta t)\, .
\end{aligned}
\end{equation}
Here we used the fact that
$\lim_{n\rightarrow\infty} \frac{1}{n}\sum_{i=0}^{n-1} f(x_i) \Delta x_i^k = O(\Delta t^{\lceil \frac{k}{2} \rceil})$,
where $\lceil \cdot \rceil$ denotes the ceiling function; this last relation  is easily verified by substituting
$\Delta x_i$ by (\ref{multi:sde:explMilstein}) and then applying the ergodic average lemmas in Appendix A.
\end{proof}

\subsubsection{Quadratic potential on $\mathbb R$}
\begin{figure}[!htb]
\begin{center}
\includegraphics[width=.8\textwidth]{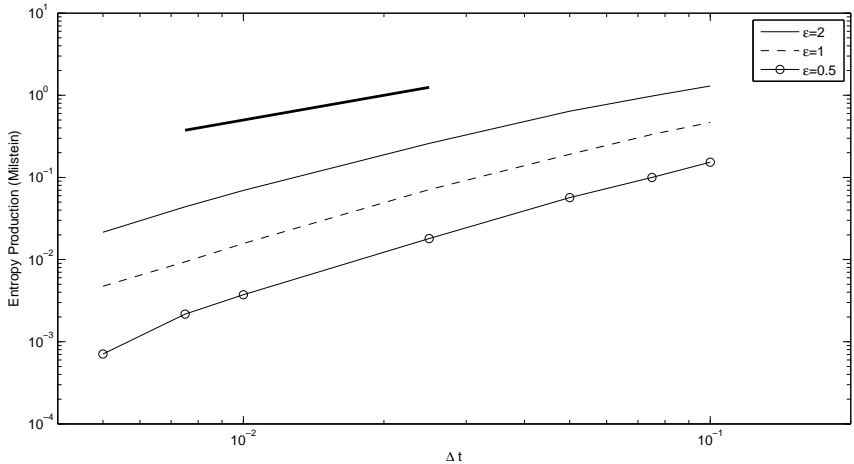}
\caption{Entropy production rate as a function of time step $\Delta t$ for the explicit Milstein's
scheme. The decrease of the entropy production rate for this numerical scheme is linear.
%Thus, in a loose sense, the reversibility property of the original continuous-time process
%is restored.
}
\label{multi:fig:Milstein}
\end{center}
\end{figure}

We compute numerically the entropy production as the time-average of the 
GC action functional. Figure~\ref{multi:fig:Milstein} shows
the numerically computed entropy production for the same example shown in
Figure~\ref{multi:fig1}. Evidently, entropy production rate decreases at least linearly
as time step $\Delta t$ is decreasing as Theorem~\ref{Milstein:entr:prod:multi:theorem}
asserts.
%Additionally, a number  of different variance functions which satisfy the condition of
%Theorem~\ref{entr:prod:multi:theorem} were tested and in all cases the decrease of
%the entropy production for the Milstein's scheme was linear. Thus, we conjecture that
%entropy production of overdamped Langevin process with multiplicative noise is of order
%$O(\Delta t)$ for Milstein's scheme.

\begin{rmrk}
We note that the rigorous asymptotics for the entropy production quickly 
become quite involved as the Milstein scheme analysis demonstrates.  
However, the GC functional is easily accessible numerically and this
allows to assess the reversibility of
each scheme computationally,  as demonstrated in  Figure~\ref{multi:fig1}
and Figure~\ref{multi:fig:Milstein}.
\end{rmrk}

\section{Entropy Production for Langevin Processes}
\label{Langevin:sec}
Let us consider another important class of reversible processes, namely the processes
driven by the Langevin equation 
\begin{equation}
\begin{aligned}
&dq_t = M^{-1}p_t dt \\
&dp_t = -\nabla V(q_t)dt - \gamma M^{-1}p_t dt + \sigma dB_t 
%\\
%&(q_0, p_0) \sim \text{invariant measure}
\label{Langevin:sde:cont}
\end{aligned}
\end{equation}
where $q_t\in\mathbb R^{dN}$ is the position vector of the $N$ particles,
$p_t\in\mathbb R^{dN}$ is the momentum vector of the particles, $M$ is
the mass matrix, $V$ is the potential energy, $\gamma$ is the friction factor (matrix),
$\sigma$ is the diffusion factor (matrix) and $B_t$ is a $dN$-dimensional Brownian
motion. Even though the Langevin system is degenerate since the noise applies
only to the momenta, the process is hypoelliptic  and is ergodic under mild conditions on $V$. 
The fluctuation-dissipation theorem asserts that friction and diffusion
terms are related with the inverse temperature $\beta\in\mathbb R$ of
the system by
\begin{equation}
(\sigma\sigma^T) = 2\beta^{-1}\gamma\, .
\end{equation}
The Langevin system is reversible (modulo momenta flip, see \eqref{DBLangevin})
with invariant measure
\begin{equation}
\mu(dq,dp) = \frac{1}{Z} \exp\left(-\beta H(q,p)\right) dq dp.
\end{equation}
where $H(q,p)$ is the Hamiltonian of the system given by 
\begin{equation}
H(q,p) = V(q) + \frac{1}{2}p^TM^{-1}p \, .
\end{equation}
Indeed if  $\mathcal L$ denotes the generator of (\ref{Langevin:sde:cont}), 
it is straightforward to verify the following modified DB condition
\begin{equation}\label{DBLangevin}
<\mathcal L f(q,p), g(q,p)>_{L^2(\mu)} =  < f(q,-p), \mathcal L g(q,-p)>_{L^2(\mu)}
\end{equation}
for any test functions $f$ and $g$ which are bounded, twice differentiable with bounded 
derivatives.  This shows that the Langevin process is reversible modulo flipping the momenta
of all particles.  

The BBK integrator \cite{Brunger:84, Lelievre:10} which utilizes a Strang splitting is applied for
the discretization of (\ref{Langevin:sde:cont}). It is written as
\begin{equation}
\begin{aligned}
p_{i+\frac{1}{2}} &= p_i -\nabla V(q_i)\frac{\Delta t}{2} - \gamma M^{-1}p_i \frac{\Delta t}{2} + \sigma \Delta W_i \\
q_{i+1} &= q_i + M^{-1}p_{i+\frac{1}{2}} \Delta t \\
p_{i+1} &= p_{i+\frac{1}{2}} -\nabla V(q_{i+1})\frac{\Delta t}{2} - \gamma M^{-1}p_{i+1} \frac{\Delta t}{2} + \sigma \Delta W_{i+\frac{1}{2}} 
%\\
%&(q_0, p_0) \sim \text{invariant measure}
\label{Langevin:BBK:integrator1}
\end{aligned}
\end{equation}
with $\Delta W_i, \Delta W_{i+\frac{1}{2}} \sim N(0, \frac{\Delta t}{2} I_{dN})$.
Its stability and convergence properties were studied in \cite{Brunger:84,
Lelievre:10} while its ergodic properties can be found in \cite{Talay:02, Mattingly:02, Mattingly:10}.
An important property of this numerical scheme which simplifies the computation of the transition 
probabilities is that the transition probabilities are non-degenerate. We rewrite the BBK integrator as 
\begin{subequations}
\begin{equation}
q_{i+1} = q_i + M^{-1}[p_i -\nabla V(q_i)\frac{\Delta t}{2} - \gamma M^{-1}p_i \frac{\Delta t}{2}] \Delta t + M^{-1} \sigma  \Delta t \Delta W_i \\
\end{equation}
\begin{equation}
p_{i+1} = (I+\gamma M^{-1}\frac{\Delta t}{2})^{-1}[\frac{1}{\Delta t}M(q_{i+1}-q_i) - \nabla V(q_{i+1})\frac{\Delta t}{2}]
+ (I+\gamma M^{-1}\frac{\Delta t}{2})^{-1} \sigma \Delta W_{i+\frac{1}{2}}
\end{equation}
\label{Langevin:BBK:integrator2}
\end{subequations}
and thus the transition probabilities of the discrete-time approximation process are 
given by the product
\begin{equation}
\Pi(q_i,p_i,q_{i+1},p_{i+1}) = P(q_{i+1}|q_i,p_i) P(p_{i+1}|q_{i+1},q_i,p_i)
\end{equation}
where $P(q_{i+1}|q_i,p_i)$ is the propagator of the positions given by
\begin{equation}
\begin{aligned}
P(q_{i+1}|q_i,p_i) &= \frac{1}{Z_0} \exp\{-\frac{1}{\Delta t^3} (\Delta q_i + M^{-1}(p_i-\nabla V(q_i)\frac{\Delta t}{2} + \gamma M^{-1}p_i\frac{\Delta t}{2})\Delta t)^T \\
&(\sigma M^{-T}M^{-1}\sigma^T)^{-1}  (\Delta q_i + M^{-1}(p_i-\nabla V(q_i)\frac{\Delta t}{2} + \gamma M^{-1}p_i\frac{\Delta t}{2})\Delta t) \} \\
\label{trans:prob:position}
\end{aligned}
\end{equation}
where $\Delta q_i = q_{i+1}-q_i$ while $P(p_{i+1}|q_{i+1},q_i,p_i)$ is the propagator of
the momenta given by
\begin{equation}
\begin{aligned}
&P(p_{i+1}|q_{i+1},q_i,p_i) = \frac{1}{Z_1} \exp\{-\frac{1}{\Delta t}
(p_{i+1} - (I+\gamma)M^{-1}\frac{\Delta t}{2})^{-1}(\frac{1}{\Delta t}M\Delta q_i - \nabla V(q_{i+1})\frac{\Delta t}{2}) )^T \\
&(\sigma^T(I+\gamma M)^{-T}(I+\gamma M^{-1})\sigma)^{-1}
(p_{i+1} - (I+\gamma)M^{-1}\frac{\Delta t}{2})^{-1}(\frac{1}{\Delta t}M\Delta q_i - \nabla V(q_{i+1})\frac{\Delta t}{2})) \} \\
\label{trans:prob:momenta}
\end{aligned}
\end{equation}
Finally, since the Langevin process is reversible modulo flip of the momenta,
the GC action functional takes the form
\begin{equation}
W(n;\Delta t) = \sum_{i=0}^{n-1} \log \frac{\Pi(q_i,p_i,q_{i+1},p_{i+1})}{\Pi(q_{i+1},-p_{i+1},q_i,-p_i)}.
\end{equation}

\subsection{Langevin Process with Additive Noise}
In the following we assume for simplicity that particles have equal masses (i.e. $M = m I$)
and that  $\sigma = \sigma I$, $\gamma = \gamma I$.  
%, even though the general case can be handled, we restrict for clarity to
%the simpler additive noise case. Thus, we assume that $\sigma(q_i) = \sigma I$,
%$\gamma(q_i) = \gamma I$ as well that particles have equal masses ($M = m I$).
In the next lemma we compute the GC action functional.

\begin{lmm}
The GC action functional of the BBK integrator equals to
\begin{equation}
W(n;\Delta t) \dot{=} \frac{\beta}{\Delta t} \sum_{i=0}^{n-1} \left[\Delta p_i^T\Delta q_i - \frac{\Delta t^2}{2m}(\nabla V(q_i)^Tp_i + \nabla V(q_{i+1})^Tp_{i+1})\right]
\label{GC:action:func:Langevin}
\end{equation}
\label{GC:action:func:Langevin:lemma}
\end{lmm}
\begin{proof}
Firstly, (\ref{trans:prob:position}) and (\ref{trans:prob:momenta}) are rewritten as
\begin{equation}
P(q_{i+1}|q_i,p_i) = \frac{1}{Z_0} \exp\left\{-\frac{m^2}{\sigma^2\Delta t^3}
|\Delta q_i + (p_i-\frac{1}{m}\nabla V(q_i)\frac{\Delta t}{2} + \frac{\gamma}{m}p_i\frac{\Delta t}{2})\Delta t|^2\right\}
\end{equation}
and
\begin{equation}
P(p_{i+1}|q_{i+1},q_i,p_i) = \frac{1}{Z_1} \exp\left\{-\frac{1}{\sigma^2\Delta t}
|(1+\frac{\gamma\Delta t}{2m})p_{i+1}-(\frac{m}{\Delta t}\Delta q_i - \frac{\Delta t}{2}\nabla V(q_{i+1})) |^2\right\}
\end{equation}
respectively.
Then, as in the overdamped Langevin case, the computation of the GC action functional is straightforward,
\begin{equation*}
\begin{aligned}
W(n; \Delta t) &= -\frac{m^2}{\sigma^2\Delta t^3} \sum_{i=0}^{n-1}
\left[\left|\Delta q_i + \frac{\Delta t^2}{2m}\nabla V(q_i) - \frac{\Delta t}{m}(1-\frac{\gamma\Delta t}{2m})p_i\right|^2
- \left|-\Delta q_i + \frac{\Delta t^2}{2m}\nabla V(q_{i+1}) + \frac{\Delta t}{m}(1-\frac{\gamma\Delta t}{2m})p_{i+1}\right|^2 \right] \\
& -\frac{1}{\sigma^2\Delta t} \sum_{i=0}^{n-1}
\left[\left|(1+\frac{\gamma\Delta t}{2m})p_{i+1}-\frac{m}{\Delta t}\Delta q_i + \frac{\Delta t}{2}\nabla V(q_{i+1})\right|^2
- \left|-(1+\frac{\gamma\Delta t}{2m})p_i+\frac{m}{\Delta t}\Delta q_i + \frac{\Delta t}{2}\nabla V(q_i)\right|^2\right] \\
&= -\frac{m^2}{\sigma^2\Delta t^3} \sum_{i=0}^{n-1} \left[
|\Delta q_i|^2 + |\frac{\Delta t^2}{2m}\nabla V(q_i)|^2 + |\frac{\Delta t}{m}(1-\frac{\gamma\Delta t}{2m})p_i|^2
+\frac{\Delta t^2}{m}\Delta q_i^T\nabla V(q_i) \right. \\
&- \frac{2\Delta t}{m}(1-\frac{\gamma\Delta t}{2m})\Delta q_i^Tp_i - \frac{\Delta t^3}{m^2}(1-\frac{\gamma\Delta t}{2m})\nabla V(q_i)^Tp_i \\
&-|\Delta q_i|^2 - |\frac{\Delta t^2}{2m}\nabla V(q_{i+1})|^2 - |\frac{\Delta t}{m}(1-\frac{\gamma\Delta t}{2m})p_{i+1}|^2
+\frac{\Delta t^2}{m}\Delta q_i^T\nabla V(q_{i+1}) \\
&\left.+ \frac{2\Delta t}{m}(1-\frac{\gamma\Delta t}{2m})\Delta q_i^Tp_{i+1} - \frac{\Delta t^3}{m^2}(1-\frac{\gamma\Delta t}{2m})\nabla V(q_{i+1})^Tp_{i+1}\right] \\
& -\frac{1}{\sigma^2\Delta t} \sum_{i=0}^{n-1}\left[
|(1+\frac{\gamma\Delta t}{2m})p_{i+1}|^2 + |\frac{m}{\Delta t}\Delta q_i|^2 + |\frac{\Delta t}{2}\nabla V(q_{i+1})|^2
- (1+\frac{\gamma\Delta t}{2m})\frac{2m}{\Delta t}p_{i+1}^T\Delta q_i \right. \\
&+ (1+\frac{\gamma\Delta t}{2m})\Delta t p_{i+1}^T\nabla V(q_{i+1}) - m\Delta q_i^T\nabla V(q_{i+1}) \\
&- |(1+\frac{\gamma\Delta t}{2m})p_i|^2 - |\frac{m}{\Delta t}\Delta q_i|^2 - |\frac{\Delta t}{2}\nabla V(q_i)|^2
+ (1+\frac{\gamma\Delta t}{2m})\frac{2m}{\Delta t}p_i^T\Delta q_i \\
&+\left. (1+\frac{\gamma\Delta t}{2m})\Delta t p_i^T\nabla V(q_i) - m\Delta q_i^T\nabla V(q_i) \right] \, .
\end{aligned}
\end{equation*}
Thus we have,
\begin{equation*}
\begin{aligned}
W(n; \Delta t)
&\dot{=} -\frac{m^2}{\sigma^2\Delta t^3} \sum_{i=0}^{n-1} \left[ \frac{\Delta t^2}{m}\Delta q_i^T(\nabla V(q_i)+\nabla V(q_{i+1}))
+\frac{2\Delta t}{m}(1-\frac{\gamma\Delta t}{2m})\Delta q_i^T\Delta p_i \right. \\
&\left.- \frac{\Delta t^3}{m^2}(1-\frac{\gamma\Delta t}{2m}) (\nabla V(q_{i+1})^Tp_{i+1}+\nabla V(q_i)^Tp_i) \right] \\
&-\frac{1}{\sigma^2\Delta t} \sum_{i=0}^{n-1}\left[-(1+\frac{\gamma\Delta t}{2m})\frac{2m}{\Delta t}\Delta p_i^T\Delta q_i
-m\Delta q_i^T(\nabla V(q_i) + \nabla V(q_{i+1})) \right. \\
&\left.+ (1+\frac{\gamma\Delta t}{2m})\Delta t (p_i^T\nabla V(q_i)+p_{i+1}^T\nabla V(q_{i+1})) \right] \\
&= \frac{2m}{\sigma^2\Delta t^2} \sum_{i=0}^{n-1} \left[ -(1-\frac{\gamma\Delta t}{2m})\Delta q_i^T\Delta p_i
+(1+\frac{\gamma\Delta t}{2m})\Delta q_i^T\Delta p_i \right. \\
&\left.+ \frac{\Delta t^2}{2m} (1-\frac{\gamma\Delta t}{2m}) (\nabla V(q_{i+1})^Tp_{i+1}+\nabla V(q_i)^Tp_i)
-\frac{\Delta t^2}{2m} (1+\frac{\gamma\Delta t}{2m}) (\nabla V(q_{i+1})^Tp_{i+1}+\nabla V(q_i)^Tp_i)
\right] \\
&= \frac{2\gamma}{\sigma^2\Delta t} \sum_{i=0}^{n-1} \left[\Delta p_i^T\Delta q_i - \frac{\Delta t^2}{2m}(\nabla V(q_{i+1})^Tp_{i+1}+\nabla V(q_i)^Tp_i)\right]
\end{aligned}
\end{equation*}
which is equal with (\ref{GC:action:func:Langevin}).
\end{proof}

\begin{rmrk} Proceeding as in Remark~\ref{over:Langevin:remark} we can compare 
the GC action functional of the BBK integrator to the GC functional for the  
additive Langevin process with constant temperature, which is given,  \cite{Maes:00}, by
\begin{equation}
W_{cont}(t) = \frac{\beta}{m} \int_0^t \nabla V(q_t)p_tdt \approx \frac{\beta\Delta t}{2m} \sum_{i=0}^{n-1} (\nabla V(q_{i+1})^Tp_{i+1} + \nabla V(q_i)^Tp_i)
\end{equation}
and is a boundary term in continuous time. Comparing the GC functionals, it is evident that the
discrete version of $W_{cont}(t)$ is contained in the functional 
$W(n;\Delta t)$ given by \eqref{GC:action:func:Langevin}.  This is similar to the overdamped 
Langevin case when discretized utilizing the explicit EM scheme. In addition the remaining term in 
the GC action functional $W(n;\Delta t)$ stems from the Strang splitting of the numerical scheme.
Moreover, this additional term critically affects the irreversibility of the discrete-time approximation process 
since it is the leading order term in the entropy production rate, as shown in the following theorem.
\end{rmrk}

\begin{thrm}
Let Assumption~\ref{basic:assumption} hold. Assume also that the potential
function $V$ has bounded fifth-order derivative.
Then, for sufficiently small $\Delta t$, there exists $C=C(N,\gamma,m)>0$ such that
\begin{equation}
EP(\Delta t) \leq C\Delta t
\label{entr:prod:add:Langevin}
\end{equation}
\label{add:noise:theorem:Langevin}
\end{thrm}

\begin{proof}
Solving (\ref{Langevin:BBK:integrator2}a) for $p_i$ and multiplying with the transpose of $p_i$,
the square of the absolute of the momenta equal to
\begin{equation}
(1-\frac{\gamma\Delta t}{2m})|p_i|^2 = \frac{m}{\Delta t}p_i^T\Delta q_i + \frac{\Delta t}{2}p_i^T\nabla V(q_i) - \sigma p_i^T\Delta W_i
\end{equation}
and similarly for $p_{i+1}$ in (\ref{Langevin:BBK:integrator2}b)
\begin{equation}
(1+\frac{\gamma\Delta t}{2m})|p_{i+1}|^2 = \frac{m}{\Delta t}p_{i+1}^T\Delta q_i + \frac{\Delta t}{2}p_{i+1}^T\nabla V(q_{i+1}) + \sigma p_{i+1}^T\Delta W_{i+\frac{1}{2}} \ .
\end{equation}
Taking the difference between the above two equations for the momenta, we obtain
\begin{equation}
|p_{i+1}|^2-|p_i|^2 + \frac{\gamma\Delta t}{2m} (|p_{i+1}|^2+|p_i|^2) = \frac{m}{\Delta t}\Delta p_i^T\Delta q_i
- \frac{\Delta t}{2}(p_{i+1}^T\nabla V(q_{i+1})+p_i^T\nabla V(q_i)) + \sigma (p_{i+1}^T\Delta W_{i+\frac{1}{2}}+ p_i^T\Delta W_i)\ ,
\end{equation}
hence, the GC action functional is rewritten as
\begin{equation}
\begin{aligned}
W(n;\Delta t) &\dot{=} \frac{\beta}{\Delta t} \sum_{i=0}^{n-1} \left[\Delta p_i^T\Delta q_i - \frac{\Delta t^2}{2m}(\nabla V(q_{i+1})^Tp_{i+1}+\nabla V(q_i)^Tp_i)\right] \\
&= \frac{\beta}{m} \sum_{i=0}^{n-1} \left[|p_{i+1}|^2-|p_i|^2 + \frac{\gamma\Delta t}{2m} (|p_{i+1}|^2+|p_i|^2)
-  \sigma (p_{i+1}^T\Delta W_{i+\frac{1}{2}}+ p_i^T\Delta W_i)  \right] \\
&\dot{=} \frac{\beta\gamma\Delta t}{m^2} \sum_{i=0}^{n-1} |p_i|^2 - \frac{\beta\sigma}{m} \sum_{i=0}^{n-1} (p_{i+1}^T\Delta W_{i+\frac{1}{2}} + p_i^T\Delta W_i)
\label{Langevin:GC:act:func2}
\end{aligned}
\end{equation}

Using the fact that $p_i$ and $\Delta W_i$ are independent while $p_{i+1}$ and $\Delta W_{i+\frac{1}{2}}$ are not
as well as the fact that the momenta in the continuous setting are zero-mean Gaussian r.v. with variance $\frac{m}{\beta}I_{dN}$,
the entropy production rate for the BBK integrator becomes
\begin{equation}
\begin{aligned}
EP(\Delta t) &= \frac{\beta\gamma}{m^2}\lim_{n\rightarrow\infty}\frac{1}{n} \sum_{i=0}^{n-1} |p_i|^2
- \frac{\beta\sigma}{m\Delta t} \lim_{n\rightarrow\infty}\frac{1}{n} \sum_{i=0}^{n-1} (p_{i+1}^T\Delta W_{i+\frac{1}{2}} + p_i^T\Delta W_i) \\
&= \frac{\beta\gamma}{m^2} \mathbb E_{\bar{\mu}}[|p|^2] -  \frac{\beta\sigma^2}{m\Delta t\big(1+\frac{\gamma\Delta t}{2m}\big)} \mathbb E_{\rho}[|\Delta W|^2] \\
&= \frac{\beta\gamma}{m^2} \left(\frac{m N}{\beta} + O(\Delta t)\right) - \frac{4\gamma}{\Delta t(2m+\gamma\Delta t)} \frac{\Delta t N}{2} \\
&= \frac{\gamma^2 N}{m(2m+\gamma\Delta t)} \Delta t + O(\Delta t)
\label{entr:prod:Langevin:final}
\end{aligned}
\end{equation}
which completes the proof.

\end{proof}

\subsubsection{Quadratic potential on a torus}
\begin{figure}[!htb]
\begin{center}
\includegraphics[width=.8\textwidth]{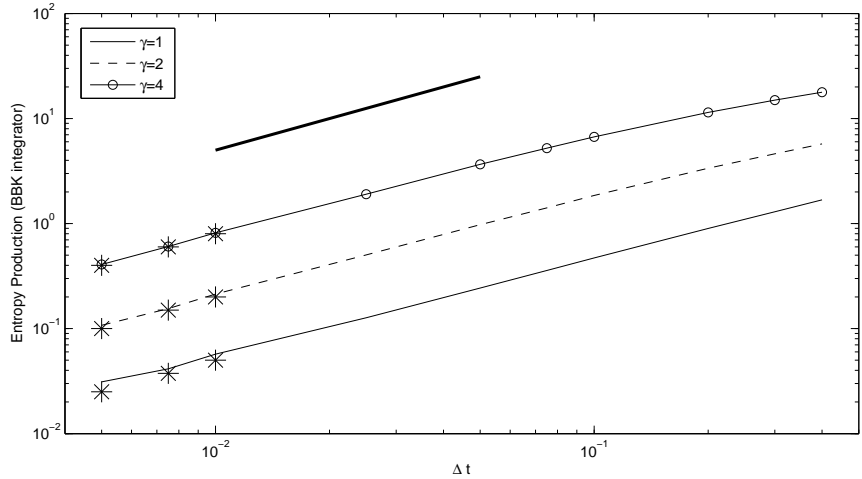}
\caption{Entropy production rate as a function of time step, $\Delta t$, for various friction
factors $\gamma$. The decrease of the entropy production rate is linear as
Theorem~\ref{add:noise:theorem:Langevin} asserts. Additionally, the theoretically-computed
entropy production rate (star points) perfectly matches the numerically-computed entropy rate.}
\label{Langevin:fig}
\end{center}
\end{figure}

The conclusions of the above theorem are illustrated by a numerical example
where the potential function is quadratic, $V(x) = \frac{|x|^2}{2}$.
Figure~\ref{Langevin:fig} shows the behavior of numerical entropy production rate
as a function of $\Delta t$ computed as the time-average of the GC action functional.
Number of particles was set to $N=5$ while the mass of its particle was set to $m=1$.
The variance of the stochastic term was set $\sigma^2=0.01$ while the final time was
set to $t=2\cdot10^5$. The initial data was chosen randomly from the zero-mean
Gaussian distribution with appropriate variance. Notice also that due to the quadratic
potential of this example Gaussian distribution is also the invariant measure of
the process. Thus, the simulation is performed at the equilibrium regime. Evidently,
the entropy production rate is of order $O(\Delta t)$ as it is expected.
Additionally, we plot (stars in the Figure) the leading term of the theoretical
value of the entropy production rate as it given by (\ref{entr:prod:Langevin:final}).
Apparently, the theoretical coefficient, $\frac{N\gamma^2}{2m^2}$, is very close to the numerically-computed
coefficient. Finally, notice that the entropy production rate is quadratically proportional
to the friction factor $\gamma$ which is in accordance with (\ref{entr:prod:Langevin:final}).

%
%\subsection{Metropolization}
%A well-known approach to make the discrete-time process which does not satisfy detailed
%balance again reversible with the desired invariant measure is to Metropolize the update
%step. Metropolized A... Langevin algorithm (MALA) was introduced by Twentie and ...
%\cite{Twentie:86} and it has been extensively studied \cite{BouRabee:09}, \cite{Lelievre:10}.
%We briefly review MALA for the overdamped Langevin process, i.e., (\ref{add:sde:cont}),
%as a solution of among others fixing the reversibility of the discrete process.
%
%In MALA, a new value, $\tilde{x}_{i+1}$, is proposed based on (\ref{add:sde:explEuler})
%and the rate
%\begin{equation}
%r(x_i,\tilde{x}_{i+1}) = e^{-\beta(V(\tilde{x}_{i+1})-V(x_i))}\frac{\Pi(\tilde{x}_{i+1},x_i)}{\Pi(x_i,\tilde{x}_{i+1})}
%\label{MALA:rate}
%\end{equation}
%is defined. Then, the proposal, $\tilde{x}_{i+1}$, is accepted with probability
%$P_i = \min(1,r(x_i,\tilde{x}_{i+1}))$ while it is rejected with probability $1-P_i$.
%Thus, the transition probability of MALA from $x_i$ to $x_{i+1}$ is
%\begin{equation}
%\Pi_{MALA}(x_i,x_{i+1}) = \Pi(x_i,x_{i+1}) P_i + (1-\int P_i \Pi(x_i,y)dy)\delta(x_{i+1}-x_i)
%\end{equation}
%which is easily shown using the identity $\min(1,r)=r^{-1}\min(1,r^{-1})$ that it satisfies
%detailed balance with invariant measure given by (\ref{inv:measure:additive}). Thus,
%MALA provides a mean not only to sample from the correct invariant measure but also
%to fix the reversibility property which was lost due to the time discretization.

\section{Summary and Future Work}
\label{conclusions}
In this paper we use the entropy production rate as a novel tool to assess quantitatively 
the (lack of)  reversibility of discretization schemes for various reversible SDE's. 
Reversibility of the discrete-time approximation process is a desirable feature when 
equilibrium simulations are performed. The entropy production rate which is defined as the 
time-average of the relative entropy between the path measure of the forward process
and the path measure of the time-reversed process is  zero when the process is reversible and positive 
when it is irreversible. Thus, it provides a way  to quantify the (ir)reversibility of the approximation process. 
Moreover, under an ergodicity assumption, the entropy production rate can be computed numerically 
on-the-fly utilizing the GC action functional. This is another attractive feature of the entropy production rate.

We computed the entropy production rate for overdamped Langevin processes both 
analytically and numerically when discretized with the explicit Euler-Maruyama scheme. One of the 
main finding in this paper  is that depending on the type of the noise --additive vs multiplicative-- the 
entropy production for the explicit EM scheme had totally different  behavior. Indeed, for additive 
noise entropy production rate is of order $O(\Delta t^2)$ while for multiplicative noise it is of 
order $O(1)$. Hence, reversibility of the discrete-time approximation process does not depend only 
on the numerical scheme but also on the intrinsic characteristics of the SDE. For the Milstein's scheme
the entropy production rate  $O(\Delta t)$  for multiplicative noise. 
Furthermore, we  computed the entropy production
rate both analytically and numerically for discretization schemes of the Langevin process with 
additive noise. Specifically, we computed the entropy production rate for the BBK integrator of the 
Langevin equation which is a quasi-symplectic splitting numerical scheme. The 
rate of entropy production was shown to be of order $O(\Delta t)$.

This paper offers a new conceptual tool for the evaluation of discretization schemes
of SDE systems simulated at the equilibrium regime. We consider only the simplest schemes here
and we will analyze in future work the behavior of  the entropy production for other 
numerical schemes  such as fully implicit EM, drift-implicit EM, higher-order schemes as well as
different kind of splitting methods. Moreover, other reversible  or even non-reversible 
processes can be analyzed in the same way, in particular extended, spatially-distributed processes.
A particularly interesting example, where the reversibility of the original system is destroyed
by numerical schemes in the form of spatio-temporal fractional step approximations of the
generator,  arises in the (partly asynchronous) parallelization of Kinetic Monte Carlo algorithms
\cite{ShimAmar05b}, \cite{AKPTX}.
Finally, another possible extension of this work is to develop adaptive schemes based on
the {\em a posteriori} simulation of  entropy production rate, which  should guarantee the reversibility
or the  approximate reversibility of the discrete-time approximation process. In this direction, 
the decomposition of entropy production functional for Metropolis-adjusted Langevin algorithms (MALA)
\cite{Roberts:96,Lelievre:10}  should be further studied and understood.

\bibliographystyle{plain}
%\bibliography{../../biblio/biblio} % you might need to change the path of the bibliography file
\bibliography{biblio}

\appendix

\section{Tools for proving Theorem~\ref{add:noise:theorem}}
\label{app:a}
%First, a generalization of the trapezoidal rule is stated and proved.
\begin{lmm}[Generalized Trapezoidal Rule]
For $k$ odd,
\begin{equation}
\begin{aligned}
V(x_{i+1})-V(x_i) &= \sum_{|\alpha|=1,3,...}^{k} C_\alpha [D^\alpha V(x_{i+1})+D^\alpha V(x_i)] \Delta x_i^\alpha \\
&+ \sum_{|\alpha|=1,3,...}^{k+2} \sum_{|\beta|=k+2-|\alpha|} B_\beta [R_\alpha^\beta(x_i,x_{i+1}) + R_\alpha^\beta(x_{i+1},x_i)]\Delta x_i^{\alpha+\beta}
\label{gen:trapez:rule}
\end{aligned}
\end{equation}
where $\alpha=(\alpha_1,...,\alpha_d)$ is a typical $d$-dimensional multi-index vector,
$D^\alpha V(x) = \frac{\partial^{|\alpha|} V}{\partial x_1^{\alpha_1}...\partial x_d^{\alpha_d}}(x)$
is the $\alpha$-th partial derivative while $x^\alpha = x_{1}^{\alpha_1}...x_{d}^{\alpha_d}$.
The coefficients $C_\alpha$ are defined recursively by
\begin{equation}
\begin{aligned}
&C_\alpha = \frac{1}{2} \ \ \ \ \ \ \ \text{for} \ \ \ |\alpha|=1 \\
&C_\alpha = \frac{1}{2}\left(\frac{1}{\alpha!} - \sum_{|\gamma|=1,3,...}^{|\alpha|-2} \frac{1}{(\alpha-\gamma)!} C_\gamma\right) \ \ \ \ \ \ \ \text{for} \ \ \ |\alpha|=3,5,...,k
\end{aligned}
\end{equation}
while the coefficients $B_\beta$ are also recursively defined by
\begin{equation}
\begin{aligned}
&B_\beta = \frac{1}{2} \ \ \ \ \text{for} \ \ \ |\beta|=0 \\
&B_\beta = -\frac{1}{2} \sum_{|\gamma|=2,4,...}^{|\beta|} \frac{1}{\gamma!} B_{\beta-\gamma} \ \ \ \ \text{for} \ \ \ |\beta|=2,4,...,k+1
\end{aligned}
\end{equation}
Finally, the remainder terms are given by \\
$R_\alpha^\beta(x_i,x_{i+1}) = \frac{|\alpha|}{\alpha!} \int_0^1 (1-t)^{|\alpha|-1}D^{\alpha+\beta}V((1-t)x_i+t x_{i+1})dt$.
\end{lmm}

\begin{proof}
The starting point is the usual Taylor series expansion around $x_i$
\begin{equation}
V(x_{i+1})-V(x_i) = \sum_{|\alpha|=1}^{k+1} \frac{1}{\alpha!} D^\alpha V(x_i)\Delta x_i^\alpha
+ \sum_{|\alpha|=k+2} R_\alpha^0(x_i,x_{i+1})\Delta x_i^\alpha
\end{equation}
and around $x_{i+1}$
\begin{equation}
V(x_{i+1})-V(x_i) = -\sum_{|\alpha|=1}^{k+1} \frac{1}{\alpha!} D^\alpha V(x_{i+1}) (-\Delta x_i)^\alpha
- \sum_{|\alpha|=k+2} R_\alpha^0(x_{i+1},x_i)(-\Delta x_i)^\alpha\, .
\end{equation}
Adding the two equations we obtain the symmetrized Taylor series expansion for $V$ given by
\begin{equation}
\begin{aligned}
&V(x_{i+1})-V(x_i) = \frac{1}{2}\sum_{|\alpha|=1,3,...}^k \frac{1}{\alpha!}[D^\alpha V(x_{i+1})+D^\alpha V(x_i)]\Delta x_i^\alpha \\
&- \frac{1}{2}\sum_{|\alpha|=2,4,...}^{k+1} \frac{1}{\alpha!}[D^\alpha V(x_{i+1})-D^\alpha V(x_i)]\Delta x_i^\alpha
+ \frac{1}{2}\sum_{|\alpha|=k+2} [R_\alpha^0(x_i,x_{i+1})+R_\alpha^0(x_{i+1},x_i)]\Delta x_i^\alpha\, .
\label{symm:Taylor:series}
\end{aligned}
\end{equation}
Moreover, generalized trapezoidal formula (\ref{gen:trapez:rule}) for $D^\alpha V$ with $|\alpha|$ even is
\begin{equation}
\begin{aligned}
&D^\alpha V(x_{i+1})-D^\alpha V(x_i) = \sum_{|\gamma|=1,3,...}^{k-|\alpha|} C_\gamma
[D^{\alpha+\gamma} V(x_{i+1})+D^{\alpha+\gamma} V(x_i)] \Delta x_i^\gamma \\
&+ \sum_{|\gamma|=1,3,...}^{k+2-|\alpha|} \sum_{|\beta|=k+2-|\alpha|-|\gamma|} B_\beta
[R_\gamma^{\alpha+\beta}(x_i,x_{i+1}) + R_\gamma^{\alpha+\beta}(x_{i+1},x_i)]\Delta x_i^{\beta+\gamma}\, .
\label{gen:trapez:rule:derivative}
\end{aligned}
\end{equation}
Hence, substituting (\ref{gen:trapez:rule:derivative}) into (\ref{symm:Taylor:series}), a recursive
Taylor series expansion
\begin{equation}
\begin{aligned}
&V(x_{i+1})-V(x_i) = \frac{1}{2}\sum_{|\alpha|=1,3,...}^k \frac{1}{\alpha!}[D^\alpha V(x_{i+1})+D^\alpha V(x_i)]\Delta x_i^\alpha \\
&- \frac{1}{2}\sum_{|\alpha|=2,4,...}^{k+1} \frac{1}{\alpha!} \sum_{|\gamma|=1,3,...}^{k-|\alpha|} C_\gamma
[D^{\alpha+\gamma} V(x_{i+1})+D^{\alpha+\gamma} V(x_i)]\Delta x_i^{\alpha+\gamma} \\
&- \frac{1}{2}\sum_{|\alpha|=2,4,...}^{k+1} \frac{1}{\alpha!} \sum_{|\gamma|=1,3,...}^{k+2-|\alpha|} \sum_{|\beta|=k+2-|\alpha|-|\gamma|} B_\beta
[R_\gamma^{\alpha+\beta}(x_i,x_{i+1}) + R_\gamma^{\alpha+\beta}(x_{i+1},x_i)]\Delta x_i^{\alpha+\beta+\gamma}\\
&+ \frac{1}{2}\sum_{|\alpha|=k+2} [R_\alpha^0(x_i,x_{i+1})+R_\alpha^0(x_{i+1},x_i)]\Delta x_i^\alpha \\
&= \frac{1}{2}\sum_{|\alpha|=1,3,...}^k \frac{1}{\alpha!}[D^\alpha V(x_{i+1})+D^\alpha V(x_i)]\Delta x_i^\alpha \\
&- \frac{1}{2}\sum_{|\alpha|=3,5,...}^{k} \sum_{|\gamma|=1,3,...}^{|\alpha|-2} \frac{1}{(\alpha-\gamma)!} C_\gamma
[D^{\alpha} V(x_{i+1})+D^{\alpha} V(x_i)]\Delta x_i^{\alpha} \\
&+ \frac{1}{2}\sum_{|\alpha|=k+2} \sum_{|\beta|=k+2-|\alpha|} [R_\alpha^\beta(x_i,x_{i+1})+R_\alpha^\beta(x_{i+1},x_i)]\Delta x_i^\alpha \\
&- \frac{1}{2}\sum_{|\alpha|=1,3,...}^{k} \sum_{|\beta|=k+2-|\alpha|} \sum_{|\gamma|=2,4,...}^{|\beta|} \frac{1}{\gamma!} B_{\beta-\gamma}
[R_\alpha^{\beta}(x_i,x_{i+1}) + R_\alpha^{\beta}(x_{i+1},x_i)]\Delta x_i^{\alpha+\beta}\\
\label{recursive:Taylor:series}
\end{aligned}
\end{equation}
is obtained after rearrangements of the sums. Equating the same powers of (\ref{recursive:Taylor:series})
and (\ref{gen:trapez:rule}), the coefficients $C_\alpha$ and $B_\beta$ are obtained.

Thus far, we presented how to compute the coefficients of the generalized trapezoidal
formula. A rigorous proof of the lemma is then easily derived by induction on the order, $k$, of
(\ref{gen:trapez:rule}) and proceeding on the reverse direction of the above formulae.
\end{proof}

\begin{lmm}
Assume that the discrete-time Markov process $x_i$ driven by
\begin{equation}
\begin{aligned}
&x_{i+1} = F(x_i, \Delta W_i)
\end{aligned}
\end{equation}
where $\Delta W_i$ are i.i.d. Gaussian random variables is ergodic with invariant
measure $\bar{\mu}$. Then,
\begin{itemize}
\item[(i)] For sufficiently smooth function $h$ we have 
\begin{equation}
\lim_{n\rightarrow\infty} \frac{1}{n}\sum_{i=0}^{n-1} h(x_i,\Delta W_i) = \mathbb E_{\bar{\mu}\times\rho} [h(x,y)]\, .
\end{equation}

\item[(ii)] For sufficiently smooth functions $f$ and $g$ we have
\begin{equation}
\lim_{n\rightarrow\infty} \frac{1}{n}\sum_{i=0}^{n-1} f(x_i)g(\Delta W_i) = \mathbb E_{\bar{\mu}} [f(x)]\mathbb E_{\rho} [g(y)]\, .
\end{equation}

\item[(iii)] For sufficiently smooth functions $f$ and $g$ and for bounded $f$ holds that
\begin{equation}
\lim_{n\rightarrow\infty} \frac{1}{n}\sum_{i=0}^{n-1} f(x_i,\Delta W_i)g(\Delta W_i) =
\mathbb E_{\bar{\mu}\times\rho} [f(x,y)]\mathbb E_{\rho} [g(y)]\, ,
\end{equation}
where $\rho$ is  the Gaussian measure.
\end{itemize}
\end{lmm}

\begin{proof}
Proving (i) is based on showing that the transition density of the joint process $z_i=(x_i,\Delta W_i)$
exists and it is positive. Both are trivial since the transition density is the product of the two densities
which are both positive. Thus, irreducibility for the joint process is proved and in combination with
stationarity,  the joint process is ergodic.
Relation (ii) is a direct consequence of (i) for $h(x,y) = f(x)g(y)$.
By denoting $\bar{f}=\mathbb E_{\bar{\mu}\times\rho} [f(x,y)]$ and $\bar{g}=\mathbb E_{\rho} [g(y)]$,
(iii) is proved by applying (i),  noting that
\begin{equation}
\begin{aligned}
&\left|\frac{1}{n}\sum_{i=0}^{n-1} f(x_i,\Delta W_i)g(\Delta W_i) - \bar{f}\bar{g}\right| \\
&= \left|\frac{1}{n}\sum_{i=0}^{n-1} f(x_i,\Delta W_i)g(\Delta W_i) - \frac{1}{n}\sum_{i=0}^{n-1} f(x_i,\Delta W_i)\bar{g}
+ \frac{1}{n}\sum_{i=0}^{n-1} f(x_i,\Delta W_i)\bar{g} - \bar{f}\bar{g}\right| \\
&\leq M |\frac{1}{n}\sum_{i=0}^{n-1}g(\Delta W_i) - \bar{g}| + |\bar{g}||\frac{1}{n}\sum_{i=0}^{n-1} f(x_i,\Delta W_i) - \bar{f}|\, ,
\end{aligned}
\end{equation}
since $f$ is bounded (i.e., $|f|\leq M$). Hence, sending $n\rightarrow\infty$, (iii) is proved.
\end{proof}

\end{document}